\documentclass[10pt]{amsart}
\usepackage[leqno]{amsmath}
\usepackage{amssymb,latexsym,soul,cite,amsthm,color,enumitem,graphicx,mathtools,microtype,accents, esint, tikz, comment} 
\usepackage[colorlinks=true,
            linkcolor=blue,
            citecolor=blue,
            urlcolor=blue]{hyperref}
\definecolor{cobalt}{rgb}{0.0, 0.28, 0.67}
\usepackage[english]{babel}
\usepackage[left=2.5cm,right=2.5cm,top=2.5cm,bottom=2.5cm]{geometry}

\usetikzlibrary{arrows.meta}
\usepackage{float}

\allowdisplaybreaks
\makeatletter
\DeclareRobustCommand*{\bfseries}{%
  \not@math@alphabet\bfseries\mathbf
  \fontseries\bfdefault\selectfont
  \boldmath
}

\numberwithin{equation}{section}

\newtheorem{theorem}{Theorem}[section]
\theoremstyle{plain}
\newtheorem{lemma}[theorem]{Lemma}
\theoremstyle{plain}

\theoremstyle{plain}
\newtheorem{corollary}[theorem]{Corollary}
\newtheorem{definition}[theorem]{Definition}
\theoremstyle{definition}
\newtheorem{remark}[theorem]{Remark}

\newcommand{\R}{{\mathbb R}}
\newcommand{\eps}{\varepsilon}
\newcommand{\mm}{\mathtt{m}}

\renewcommand{\epsilon}{\varepsilon}
\renewcommand{\rho}{\varrho}
\newcommand{\beq}{\begin{equation}}
\newcommand{\eeq}{\end{equation}}
\renewcommand{\le}{\leqslant}

\renewcommand{\d}{\mathrm{d}}                 
\newcommand{\dx}{\mathrm{d}x}
\newcommand{\dy}{\mathrm{d}y}
\newcommand{\dz}{\mathrm{d}z}
\newcommand{\dt}{\mathrm{d}t}

\newcommand{\loc}{\mathrm{loc}}
\DeclareMathOperator*{\esssup}{ess\,sup}

\def\Xint#1{\mathchoice
{\XXint\displaystyle\textstyle{#1}}%
{\XXint\textstyle\scriptstyle{#1}}%
{\XXint\scriptstyle\scriptscriptstyle{#1}}%
{\XXint\scriptscriptstyle\scriptscriptstyle{#1}}%
\!\int}
\def\XXint#1#2#3{{\setbox0=\hbox{$#1{#2#3}{\int}$ }
\vcenter{\hbox{$#2#3$ }}\kern-.6\wd0}}
\def\fint{\Xint-}

\def\Yint#1{\mathchoice
    {\YYint\displaystyle\textstyle{#1}}%
    {\YYint\textstyle\scriptstyle{#1}}%
    {\YYint\scriptstyle\scriptscriptstyle{#1}}%
    {\YYint\scriptscriptstyle\scriptscriptstyle{#1}}%
      \!\iint}
\def\YYint#1#2#3{{\setbox0=\hbox{$#1{#2#3}{\iint}$}
    \vcenter{\hbox{$#2#3$}}\kern-.51\wd0}}
\def\longdash{{-}\mkern-3.5mu{-}} 

\def\fiint{\Yint\longdash}

\makeatletter
\newcommand{\leqnomode}{\tagsleft@true}
\newcommand{\reqnomode}{\tagsleft@false}
\makeatother

\usepackage{graphicx} 

\title{No Time Gap in Harnack Inequalities for Nonlocal Doubly Nonlinear Parabolic Equations}
\author[Cassanello]{Filippo Maria Cassanello}  \address{Filippo Maria Cassanello\\Università degli Studi di Cagliari, Via Ospedale 72, 09124, Cagliari, Italy }
\email{\url{filippom.cassanello@unica.it}}
\author[De Filippis]{Filomena De Filippis}  \address{Filomena De Filippis \\ Fachbereich Mathematik, Universität Salzburg, Hellbrunner Str. 34, 5020 Salzburg, Austria }
\email{\url{filomena.defilippis@plus.ac.at}}
\author[Stanko]{Calvin Stanko}  \address{Calvin Stanko\\Fachbereich Mathematik, Universität Salzburg, Hellbrunner Str. 34, 5020 Salzburg, Austria }
\email{\url{calvin.stanko@plus.ac.at}}

\begin{document}

\subjclass[2020]{\vspace{1mm} 35B65, 35R11, 76S05} 

\keywords{\vspace{1mm} Harnack inequality, doubly nonlinear, nonlocal}

\thanks{{\it Acknowledgements.} This research was funded in whole or in part by the Austrian Science Fund (FWF) [10.55776/PAT1850524]. For open access purposes, the author has applied a CC BY public copyright license to any author accepted manuscript version arising from this submission. \\ The first author is  supported by INdAM - GNAMPA Project, CUP E5324001950001. Moreover, he is also partially supported by the research project \textit{Partial Differential Equations and their Role in Understanding Natural Phenomena}
(Fondazione di Sardegna 2023, CUP F23C25000080007) \\ The second author is  supported by INdAM - GNAMPA Project, CUP E53C25002010001.}

\begin{abstract}
We establish Harnack inequalities without a time gap for nonlocal doubly nonlinear parabolic equations for every choice of $p>1$ and $0<q<\frac{N(p-1)}{(N-sp)_+}$, under optimal tail assumptions.
\end{abstract}

\maketitle
\begin{center}
\begin{minipage}{12cm}
\tableofcontents
\end{minipage}
\end{center}

\section{Introduction}
\noindent
We investigate Harnack inequalities for weak solutions to a class of
nonlocal doubly nonlinear parabolic equations. More precisely, we consider
the equation
\begin{equation}\label{eq:PDE}
    \partial_t u^q-\mathcal{L}_{K,p}^s u=0
    \qquad\mbox{in }\Omega_T,
\end{equation}
where $s\in(0,1)$, $p>1$, $q>0$, and
$\Omega_T:=\Omega\times(0,T]$ is a space-time cylinder
over an open and bounded set $\Omega\subset\mathbb{R}^N$,
with $N\geq2$. The nonlocal operator $\mathcal{L}_{K,p}^s$ is formally
defined by
\[
    -\mathcal{L}_{K,p}^s u(x,t)
    :=
    2\,\mathrm{P.V.}\int_{\mathbb{R}^N}
    |u(x,t)-u(y,t)|^{p-2}
    (u(x,t)-u(y,t))K(x,y)\,\dy,
\]
where $\mathrm{P.V.}$ denotes the principal value. The kernel
\(
    K:\mathbb{R}^N\times\mathbb{R}^N
    \longrightarrow[0,\infty)
\)
is measurable and satisfies
\begin{equation}\label{k}
    \frac{\texttt{L}^{-1}}{|x-y|^{N+sp}}
    \leq K(x,y)=K(y,x)
    \leq \frac{\texttt{L}}{|x-y|^{N+sp}},
\end{equation}
for almost every
$(x,y)\in\mathbb{R}^N\times\mathbb{R}^N$, for some constant
$\texttt{L}\geq1$.

We consider the nonlocal Trudinger equation and both the fast and slow diffusion regimes, corresponding
to $q>p-1$ and $q<p-1$, respectively. For the fast diffusion regime,
we impose the following upper restriction on $q$:
\begin{equation}\label{suq}
    q < \frac{N(p-1)}{(N-sp)_{+}},
\end{equation}
that is automatically satisfied when $q \leq p-1$.
Our main result establishes the following Harnack inequality without a time gap in all diffusion regimes. We further show that the temporal configuration of this inequality is sharp. 

\begin{theorem}[Harnack inequality without time gap]\label{timegapno}
Let $p>1$ and $q>0$ satisfy \eqref{suq}. Let $u$ be a continuous, globally nonnegative weak
solution to \eqref{eq:PDE} in $\Omega_T$ and let $(x_o,t_o)\in\Omega_T$ be such that
$u_o:=u(x_o,t_o)>0$. There exist constants $\mathfrak{s}\in(0,1)$, $\mathfrak{c}>1$ and
$\mathfrak{a}>1$, depending only
on the data, such that if
\[
B_{16\varrho}(x_o)\times
\left[
t_o-\mathfrak{a} u_o^{q+1-p}\varrho^{sp},
t_o+\mathfrak{a} u_o^{q+1-p}\varrho^{sp}
\right]\Subset\Omega_T
\]
then
\begin{equation}\label{supinf}
\sup_{B_\varrho(x_o)\times
(t_o-\mathfrak{s} u_o^{q+1-p}\varrho^{sp},t_o)}u
\leq
\mathfrak{c}
\inf_{B_\varrho(x_o)\times
(t_o,t_o+\mathfrak{s} u_o^{q+1-p}\varrho^{sp})}u.
\end{equation}
\end{theorem}
\noindent
We point out that the continuity assumed in the theorem is mainly used in a qualitative way. In particular, it makes the pointwise value \(u_o\) well defined. Note that, continuity results for related nonlocal equations have been established in \cite{ZL,NC}.

The proofs of the two regimes, fast and slow, follow different approaches. Starting with the slow diffusion regime, we first establish
intrinsic integral representation estimates, which are then used
to deduce \eqref{supinf}. Such estimates are also obtained
in the complementary range, as stated in the following result.

\begin{theorem}[Intrinsic representation estimates]\label{integralrep}
Let $p>1$ and $q>0$. Let $u$ be a continuous,
globally nonnegative weak solution to \eqref{eq:PDE} and let
$(x_o,t_o)\in\Omega_T$ be such that $u_o:=u(x_o,t_o)>0$.
There exist constants ${\mathfrak{s}}\in(0,1)$, $\mathfrak{c}>1$ and $\mathfrak{a}>1$,
depending only on the data, such that if
\[
B_{16\varrho}(x_o)\times
\left[t_o-\mathfrak{a} u_o^{q+1-p}\varrho^{sp},\,
t_o+\mathfrak{a} u_o^{q+1-p}\varrho^{sp}\right]\Subset\Omega_T,
\]
then
\begin{align}\label{integralbounds}
&\frac{1}{\mathfrak{c}}
\int_{t_o-{\mathfrak{s}} u_o^{q+1-p}\varrho^{sp}}^{t_o}
\int_{\R^N}
\frac{u^{p-1}(x,t)}{(\varrho+|x-x_o|)^{N+sp}}\,\dx\,\dt
\leq u_o^q.
\end{align}
Moreover, if $q$ is as in \eqref{suq}, it holds
\begin{align} \label{integralbounds2}
    u_o^q\leq\mathfrak{c}
\int_{t_o-{\mathfrak{s}} u_o^{q+1-p}\varrho^{sp}}^{t_o}
\int_{\R^N}
\frac{u^{p-1}(x,t)}{(\varrho+|x-x_o|)^{N+sp}}\,\dx\,\dt.
\end{align}
\end{theorem}
\noindent
 All results for the nonlocal Trudinger equation follow from the slow diffusion approach. Let us now describe the two approaches in slightly more detail.
In both diffusion regimes, the lower integral estimate \eqref{integralbounds}  is
based on a positivity estimate for the solution that follows from
Lemma~\ref{lem3.2} below. As for the upper integral estimate \eqref{integralbounds2},
in the slow diffusion regime it follows from a local
$L^\infty$-estimate, whereas in the fast diffusion regime
it is deduced from the Harnack inequality \eqref{supinf}.
Thus, in a sense, the order of the argument is reversed:
in the fast diffusion regime we use the Harnack inequality
to derive an integral estimate, while in the slow diffusion
regime and $q=p-1$ the integral estimates are used to prove the
Harnack inequality.

The proof adopted for the Harnack inequality \eqref{supinf}
in the fast diffusion regime is considerably more delicate. The strategy combines two approaches,
one of local and the other of nonlocal nature. The local approach, based
on the integral Harnack inequality and the expansion of
positivity, is particularly effective when
the positive tail
\begin{equation} \label{tailll}
    \operatorname{Tail}_{p-1}
(u_\pm;B_R(x_o) \times (t_o-S,t_o])
:=
\int_{t_o-S}^{t_o}
\int_{\mathbb R^N\setminus B_R(x_o)}
\frac{u_\pm^{p-1}(y,t)}
{|y-x_o|^{N+sp}}\,\dy\,\dt
\end{equation}
of the solution is sufficiently small. In this case, suitable modifications to adapt the local argument to the nonlocal setting leads to a time-insensitive
Harnack inequality, as in \cite{B1}, see also \cite{B2,DGV1,DGV2, DGV3}.  When the positive tail is large, a genuinely nonlocal
approach allows us to obtain the desired Harnack inequality.
This approach is strongly inspired by the work of
Liao~\cite{NN}. Its key ingredient is again the positivity
estimate established in Lemma~\ref{lem3.2}, which relies on a refined Caccioppoli inequality in
Lemma~\ref{caccC} that makes use of moving
levels. The construction of these levels reflects the
influence of the nonlocal interactions. More precisely,
the moving levels are governed by the difference of two time-dependent
quantities: the first one involves the positive tail
and the other the negative tail. Their interplay
yields a quantitative lower bound for the solution that is the key step to prove \eqref{supinf},
provided the negative contribution is sufficiently small. In the theorems stated above, this last condition  is always satisfied due to the globally nonnegativity of the solution. However,  some of our results extend to the more general setting of locally nonnegative solution, as detailed in the following remark.
\begin{remark}[Local nonnegativity and time-dependent kernels]\label{rem:kernel-time-dependence}
In the regime $0<q\leq p-1$, the global nonnegativity assumption in Theorem~\ref{timegapno}
can be replaced by local nonnegativity, provided that
\begin{equation}\label{negativetaill}
\int_{t_o-\mathfrak{a} u_o^{q+1-p}\varrho^{sp}}^{t_o+\mathfrak{a} u_o^{q+1-p}\varrho^{sp}}
\int_{\R^N\setminus B_\varrho(x_o)}
\frac{u_-^{p-1}(x,t)}{|x-x_o|^{N+sp}}\,\dx\,\dt
\leq \mathfrak{m} u_o^q,
\end{equation}
for a sufficiently small $\mathfrak{m}\in(0,1)$ depending only on the data, where $\mathfrak{a} >1$ is the constant in Theorem \ref{timegapno}.
The lower bound in Theorem~\ref{integralrep}, namely  \eqref{integralbounds}, also holds for locally nonnegative
solutions, with $u_+^{p-1}$ in the nonlocal integral, provided that
\begin{equation}\label{rep-negative-tail-condition}
\int_{t_o-\mathfrak{s}u_o^{q+1-p}\varrho^{sp}}^{t_o}
\int_{\R^N\setminus B_{4\varrho}(x_o)}
\frac{u_-^{p-1}(x,t)}{|x-x_o|^{N+sp}}\,\dx\,\dt
\leq \mathfrak{m}u_o^q,
\end{equation}
where $\mathfrak{s}\in(0,1)$ is the constant in Theorem~\ref{integralrep}.
The same conclusions remain valid for kernels $K=K(x,y,t)$ that are
measurable in time and satisfy the bounds in \eqref{k} uniformly for
almost every $t$. The time-independence of the kernel is used only in the proof of Theorem \ref{timegapno} for the fast diffusion regime,
specifically in the comparison argument based on the principle from
\cite{Str26}.
\end{remark}
\noindent
Combining \eqref{supinf} and \eqref{integralbounds}, we obtain
a full Harnack inequality without time gap, which also includes
the positive tail over the past intrinsic cylinder on the
left-hand side; see Theorem~\ref{fullharnack}.

Some of the tools used in this paper, namely local boundedness,
the integral Harnack inequality, and the expansion of positivity,
were established for the fractional porous medium equation,
corresponding to $p=2$, in the fast diffusion regime by the second author in \cite{FDF}
and by all three authors in \cite{CDS}, respectively.
These results were subsequently extended to the doubly nonlinear
setting by \cite{SSY2026,SZ26}, to which we shall refer throughout
the paper. 
Since we also consider the slow diffusion regime, we rely on some results from \cite{SSY2026_slow}, where a Harnack inequality with a time gap is established. A corresponding result in the fast diffusion regime, under the additional restriction $\Lambda_{p-1}>0$ (for the definition see \eqref{LAMBDA} below) is obtained in \cite{SSY2026}. In Corollary~\ref{corollary}, we extend these results to a Harnack inequality without a time gap, under the same additional restriction $\Lambda_{p-1}>0$ and within the framework considered in those works, namely for local weak solutions and kernels that are  measurable in time.

In the present paper, the bound on $q$ in \eqref{suq} for the fast diffusion regime corresponds precisely to
its counterpart in the local theory of Harnack inequalities, i.e.
\begin{equation}\label{onqq}
   q<\frac{N(p-1)}{(N-p)_{+}},
\end{equation}
as established in \cite{B1, DK92}. However,
the time-insensitive behaviour of the Harnack inequality in the local fast diffusion regime is lost in the nonlocal setting, as shown by the counterexample constructed in the last section. In other words, the comparisons within the same intrinsic space-time cylinder is no longer valid.  On the other hand, our result in the slow diffusion regime and for the nonlocal Trudinger equation is rather surprising when compared with the local theory, where a waiting time is necessary. This local behaviour is illustrated in \cite{DGV3}, where for the degenerate parabolic $p$-Laplacian, a solutions with compact support exhibiting a finite speed of propagation, preventing an immediate expansion of positivity, was constructed.

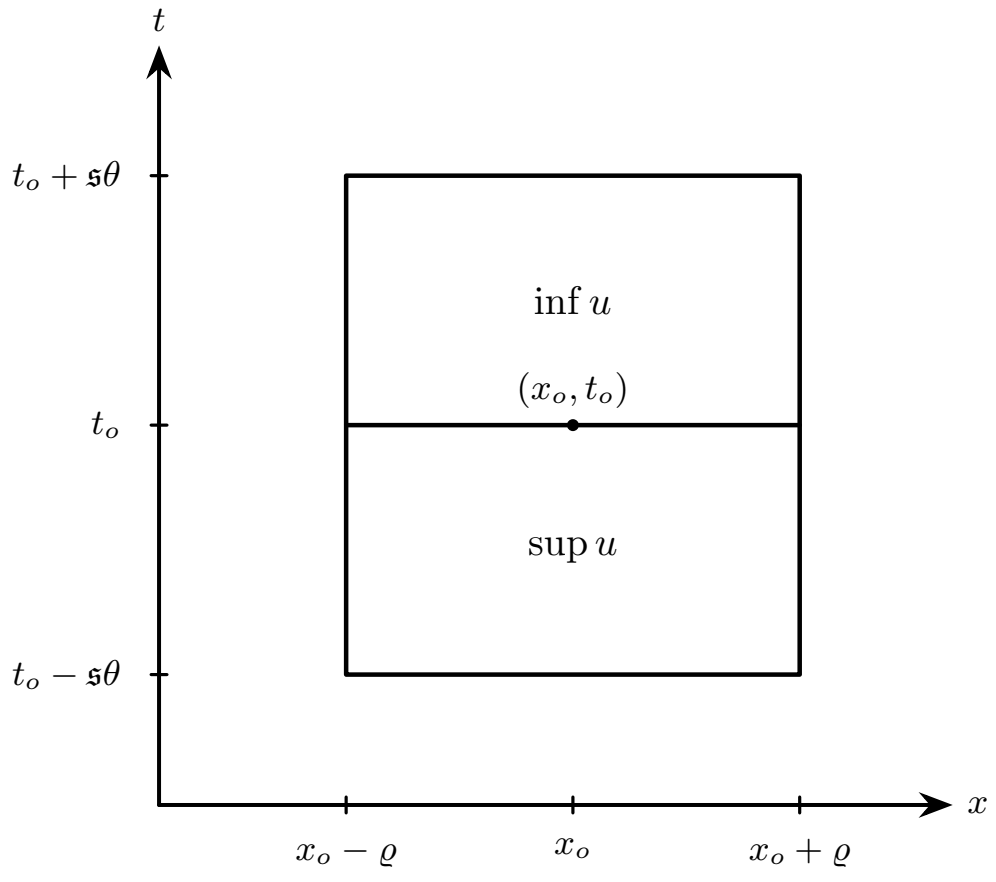
\begin{figure}[!htbp] \label{figure}
\centering
\makebox[\linewidth][c]{%
\hspace*{-5mm}%
\begin{tikzpicture}[
  scale=0.75,
  transform shape,
  >=Stealth,
  line cap=round,
  line join=round,
  every node/.style={font=\small}
]
  \def\xa{3.2} \def\xb{7.2} \def\xc{5.2}
  \def\tb{1.90} \def\tm{4.10} \def\tt{6.30}

  \draw[->, line width=.75pt] (1.55,.75) -- (8.55,.75)
    node[right] {$x$};
  \draw[->, line width=.75pt] (1.55,.75) -- (1.55,7.45)
    node[above] {$t$};

  \draw[line width=.85pt] (\xa,\tb) rectangle (\xb,\tt);
  \draw[line width=.85pt] (\xa,\tm) -- (\xb,\tm);

  \node[font=\normalsize] at (\xc,3.00) {$\sup u$};
  \node[font=\normalsize] at (\xc,5.20) {$\inf u$};

  \foreach \yy in {1.9,4.1,6.3}
    \draw[line width=.6pt] (1.48,\yy) -- (1.62,\yy);

  \node[anchor=east] at (1.34,\tb) {$t_o-\mathfrak{s}\theta$};
  \node[anchor=east] at (1.34,\tm) {$t_o$};
  \node[anchor=east] at (1.34,\tt) {$t_o+\mathfrak{s}\theta$};

  \foreach \xx in {3.2,5.2,7.2}
    \draw[line width=.6pt] (\xx,.68) -- (\xx,.82);

  \node[anchor=north] at (\xa,.58) {$x_o-\rho$};
  \node[anchor=north] at (\xc,.58) {$x_o$};
  \node[anchor=north] at (\xb,.58) {$x_o+\rho$};


  \fill (\xc,\tm) circle (1.5pt);
  \node[above, inner sep=4pt] at (\xc,\tm)
    {$(x_o,t_o)$};

\end{tikzpicture}%
\hspace*{5mm}%
}

\caption{Here
$u_o=u(x_o,t_o)>0$ and
$\theta=u_o^{\,q+1-p}\rho^{sp}$
is the intrinsic time scale.}

\end{figure}

Hence, in  the nonlocal setting, the distinction between the regimes disappears as far as the time gap is concerned: all the regimes admit a Harnack inequality without a time gap. Figure~\ref{figure} depicts the geometric configuration of this inequality: the past and future intrinsic cylinders share $t_o$ as their common temporal boundary. The supremum of $u$ over the past cylinder is controlled by the infimum over the future one. 

The example in the final section shows the sharpness of the temporal configuration of the Harnack inequality obtained in Theorem~\ref{timegapno} for all regimes. The underlying philosophy behind the construction is that a nonlocal operator detects a perturbation even at points where the function is locally constant, whereas a local differential operator does not. Consequently, this allows to consider functions that can exhibit arbitrarily large temporal oscillations locally, and so are not subject to a time-insensitive Harnack estimate, but are still solutions of the equation because of a perturbation far from the origin that compensates for the time derivative. This is a natural strategy for nonlocal operators and cannot work in the local setting, as we will show later.

\vskip 0.5 cm
\noindent
 We conclude with a brief overview of the literature on Harnack inequalities in the nonlocal setting. For linear nonlocal equations, a full parabolic Harnack inequality proved in \cite{KW} by combining local boundedness for subsolutions with a weak Harnack inequality for supersolutions. This approach was subsequently used in \cite{LW25} to establish a time-insensitive Harnack inequality. Local boundedness, Harnack inequalities and related regularity estimates can also be found in \cite{P, CCV11,FK13,Stromqvist19,Kim19,L2,Nq,CCI}. For nonlinear nonlocal parabolic equations, including fractional $p$-Laplacian and doubly nonlinear models, we refer to \cite{BGK22,StromqvistLB19,DZZ21,ByunKim24,L1,DieningKimLeeNowak25,NNA}. Further developments concern kinetic equations with nonlocal diffusion \cite{ImbertSilvestre20,AnceschiPalatucciPiccinini26,KassmannWeidnerKinetic24,KW26c} and mixed local-nonlocal parabolic equations \cite{GarainKinnunen23,Nakamura23}. Fractional porous medium equations have been investigated, in particular, in \cite{BV15,BSV15,BFV18}.
\\ \\ 
\noindent
The paper is organized as follows. Section~\ref{sec2} introduces the notation and collects the preliminary results. Section~\ref{sec5} establishes the refined Caccioppoli inequality and the positivity estimate using moving levels. In Section~\ref{sec3}, we prove the Harnack inequality without a time gap in both diffusion regimes, Theorem~\ref{timegapno}, and the intrinsic representation estimates, Theorem~\ref{integralrep}. Finally, Section~\ref{sec4} presents the counterexample showing the sharpness of the temporal configuration of the Harnack inequality.

\section{Preliminaries} \label{sec2}
\noindent
In this section, we introduce the notation and establish some preliminary results that will be used in the sequel.

\subsection*{Notation}

Throughout the paper, we denote by $\gamma$ a positive constant,
which may vary from line to line and depends only on the
structural data of the problem, namely
\[
N,s,p,q,\texttt{L}.
\]
We also use the notation $\lesssim$, $\gtrsim$, and $\approx$,
with subscripts whenever necessary, to specify the dependence
of the implicit constants. We write $\mathbb{N}_0:=\mathbb{N}\cup\{0\}$.
For $x_o\in\R^N$ and $R>0$, let $B_R(x_o)$ denote the open
ball centered at $x_o$ with radius $R$.
When the center is the origin, we simply write $B_R$.
Given $(x_o,t_o)\in\R^N\times\R$ and $R,S>0$,
we introduce the backward parabolic cylinder
\[
Q_{R,S}(x_o,t_o)
:=B_R(x_o)\times(t_o-S,t_o].
\]
For $\texttt{r}>0$, set
\begin{equation} \label{LAMBDA}
    \Lambda_{\texttt{r}}:=sp\texttt{r}+N(p-q-1).
\end{equation}
For $v\in L^1(\Omega_T,\R^k)$ and any measurable set
$U\subset\Omega$ with $|U|>0$, we denote by $[v]_U$
the spatial average of $v$ over $U$, defined for
almost every $t\in(0,T)$ by
\[
[v]_U(t)
:=\fint_U v(x,t)\,\dx
=\frac{1}{|U|}\int_U v(x,t)\,\dx.
\]
Analogously, if $Q\subset\Omega_T$ is measurable
and $|Q|>0$, the space-time average of $v$ over $Q$
is given by
\begin{equation}\label{vq}
[v]_Q
:=\fiint_Q v\,\dx\,\dt
=\frac{1}{|Q|}\iint_Q v\,\dx\,\dt.
\end{equation}
For $s\in(0,1)$ and $p>1$, we consider the fractional
Sobolev space
\[
W^{s,p}(\R^N)
:=
\left\{
v\in L^p(\R^N):
\iint_{\R^N\times\R^N}
\frac{|v(x)-v(y)|^p}{|x-y|^{N+sp}}
\,\dx\,\dy<\infty
\right\},
\]
equipped with the norm
\[
\|v\|_{W^{s,p}(\R^N)}
:=
\left(\int_{\R^N}|v|^p\,\dx\right)^{\frac1p}
+
\left(
\iint_{\R^N\times\R^N}
\frac{|v(x)-v(y)|^p}{|x-y|^{N+sp}}
\,\dx\,\dy
\right)^{\frac1p}.
\]
The space $W^{s,p}(\Omega)$ is defined in the same way,
with the integrals restricted to $\Omega$.
We further set
\[
W^{s,p}_0(\Omega)
:=
\left\{
v\in W^{s,p}(\R^N):
v=0 \quad\text{a.e. in }\R^N\setminus\Omega
\right\}.
\]
We are now ready to introduce the notion of
local weak super(sub)-solution. 
\begin{definition} \label{def}
A measurable function $u: \R^N \times (0,T] \to \R$ satisfying
$$ u \in C_{\loc} ((0,T]; L^{q+1}_{\loc}( \Omega) ) \cap L^p_{\loc} ((0,T]; W^{s,p}_{\loc}( \Omega) ), $$
is a local, {weak super(sub)-solution} to the equation  \eqref{eq:PDE} in $\Omega_T$, if for every sub-interval $[t_1,t_2] \subset (0,T]$, every bounded open set $\tilde{\Omega} \Subset \Omega$, it holds
\begin{align} \label{ineq:L_1_tail_condition}
    \int_{t_1}^{t_2} \int_{\R^N} \frac{|u(x,t)|^{p-1}}{1+|x|^{N+sp}}\,\dx\,\dt < \infty
\end{align}
and
\begin{align}  \label{eq}
& \int_{{\tilde{\Omega}}} |u|^{q-1}u \phi\,\dx \bigg |_{t_1}^{t_2}-
\int_{t_1}^{t_2}\int_{ {\tilde{\Omega}}}  |u|^{q-1}u \, \partial_t\phi\,\dx\,\dt \notag \\ & \qquad \qquad + \int_{t_1}^{t_2} \iint_{\R^N \times \R^N}|u(x,t)-u(y,t)|^{p-2}(u(x,t)-u(y,t))(\phi(x,t)-\phi(y,t)) K(x,y)\,\dx\, \dy \,\dt \geq (\leq) \, 0
\end{align}
is satisfied for all nonnegative test functions
\[
\phi \in  L^p_{\loc} ((0,T]; W^{s,p}_0({\tilde{\Omega}}) ) \cap  W^{1,q+1}_{\loc} ((0,T]; L^{q+1}({\tilde{\Omega}})) .
\]
A function that is both a weak super- and sub-solution is called a weak solution.
\end{definition}
Note that \eqref{ineq:L_1_tail_condition} ensures the finiteness of the tail defined in \eqref{tailll}. For the introduction of the notion of nonlocal tail, we refer to \cite{DKP14,DKP16}.
\subsection*{Useful results}

For $u,k\in\mathbb{R}$ and $q>0$, we introduce
\begin{align}\label{GG}
    \mathfrak{g}_{\pm}(u,k):=
    \pm q\int_k^u |\theta|^{q-1}
    (\theta-k)_{\pm}\,\d\theta,
\end{align}
where
\begin{align*}
    (\theta-k)_+:=\max\{\theta-k,0\},
    \qquad
    (\theta-k)_-:=\max\{-(\theta-k),0\}.
\end{align*}
The following elementary estimate can be found
in \cite[Lemma 2.2]{BDL_21}.
\begin{lemma}\label{2.3}
Let $q>0$ and $\mathfrak{g}_{\pm}$ be defined
as in \eqref{GG}. Then there exists a constant
$c=c(q)>0$ such that, for every $a,b\in\mathbb{R}$,
\[
    \frac{1}{c}(|a|+|b|)^{q-1}(a-b)_{\pm}^2
    \leq \mathfrak{g}_{\pm}(a,b)
    \leq c(|a|+|b|)^{q-1}(a-b)_{\pm}^2.
\]
\end{lemma}
\noindent
We also recall the fast geometric convergence
lemma, see \cite[Lemma 7.1]{G}.
\begin{lemma}\label{fg}
    Let $(Y_n)_{n\in\mathbb{N}_0}$ be a sequence
    of positive real numbers satisfying
    \[
        Y_{n+1}\leq c\,b^n\,Y_n^{1+\alpha},
    \]
    where $c,b>1$ and $\alpha>0$ are fixed.
    If
    \[
        Y_0\leq c^{-1/\alpha}b^{-1/\alpha^2},
    \]
    then $Y_n\to0$ as $n\to\infty$.
\end{lemma}
\noindent
We next state a parabolic Sobolev embedding, see \cite[Lemma 2.2]{SZ26}.
\begin{lemma}[Parabolic embedding]\label{em}
Let $s\in(0,1)$, $p>1$, $q>0$, and set
\begin{equation}\label{m}
\textnormal{\texttt{m}}
:=p+\frac{sp}{N}(q+1)
=\frac{p}{N}\left(N+s(q+1)\right).
\end{equation}
Let $B_R\subset\mathbb R^N$ be a ball and
$t_1<t_2$. Suppose that
\[
u\in L^p\left(t_1,t_2;W^{s,p}(B_R)\right)
\cap L^\infty\left(t_1,t_2;L^{q+1}(B_R)\right)
\]
and that, for some $d\in(0,1)$,
$\operatorname{supp}u(\cdot,t)\Subset B_{(1-d)R}$
for almost every $t\in(t_1,t_2)$. Then
\begin{align}\label{eq:parabolic-sobolev}
&\int_{t_1}^{t_2}\int_{B_R}
|u(x,t)|^\textnormal{\texttt{m}}\,\dx\,\dt
\notag\\
&\quad\leq c\left(
R^{sp}\int_{t_1}^{t_2}
\iint_{B_R\times B_R}
\frac{|u(x,t)-u(y,t)|^p}{|x-y|^{N+sp}}
\,\dx\,\dy\,\dt
+\frac{1}{d^{N+sp}}
\int_{t_1}^{t_2}\int_{B_R}
|u(x,t)|^p\,\dx\,\dt
\right)
\notag\\
&\qquad\quad\times
\left(
\esssup_{t_1<t<t_2}
\fint_{B_R}|u(x,t)|^{q+1}\,\dx
\right)^{\frac{sp}{N}},
\end{align}
where $c=c(N,s,p,q)>0$.
\end{lemma}
\noindent
We now recall a Caccioppoli-type energy estimate.
The following inequalities extend the Caccioppoli
estimates used in \cite{CDS} to the nonlocal
doubly nonlinear setting. Their proof follows by a straightforward adaptation of the argument in \cite{CDS} to the doubly nonlinear setting.

\begin{lemma}[Caccioppoli estimates]\label{cacc}
Let $p>1$ and $q>0$, and let $u$ be a locally nonnegative
weak sub(super)-solution to \eqref{eq:PDE} in $\Omega_T$, with kernel $K=K(x,y,t)$ measurable in the time variable and satisfying \eqref{k}.
Assume $Q_{R,S}(x_o,t_o)\Subset\Omega_T$.
Let $k:[t_o-S,t_o]\to[0,\infty)$ be an absolutely continuous function
and set $w_\pm(x,t):=(u(x,t)-k(t))_\pm$. For every time-independent, nonnegative Lipschitz
cut-off function $\xi$ compactly supported in $B_R(x_o)$
and every $\tau\in(t_o-S,t_o]$, we have
\begin{align}\label{cacc2}
&\int_{B_R(x_o)}
\xi^p(x)\mathfrak{g}_\pm(u(x,\tau),k(\tau))\,\dx
\notag\\
&\quad+\frac{c_p}{\textnormal{\texttt{L}}}
\int_{t_o-S}^{\tau}
\iint_{B_R(x_o)\times B_R(x_o)}
\frac{|\xi(x)w_\pm(x,t)-\xi(y)w_\pm(y,t)|^p}
{|x-y|^{N+sp}}\,\dx\,\dy\,\dt
\notag\\
&\quad+\frac{1}{\textnormal{\texttt{L}}}
\int_{t_o-S}^{\tau}\int_{B_R(x_o)}
\xi^p(x)w_\pm(x,t)
\left(
\int_{\R^N}\frac{w_\mp^{p-1}(y,t)}{|x-y|^{N+sp}}\,\dy
\right)\dx\,\dt
\notag\\
&\leq \int_{B_R(x_o)}
\xi^p(x)\mathfrak{g}_\pm(u(x,t_o-S),k(t_o-S))\,\dx
\notag\\
&\quad+\gamma
\int_{t_o-S}^{\tau}
\iint_{B_R(x_o)\times B_R(x_o)}
\max\{w_\pm^p(x,t),w_\pm^p(y,t)\}
\frac{|\xi(x)-\xi(y)|^p}{|x-y|^{N+sp}}
\,\dx\,\dy\,\dt
\notag\\
&\quad+4\textnormal{\texttt{L}}
\int_{t_o-S}^{\tau}\int_{B_R(x_o)}
\xi^p(x)w_\pm(x,t)
\left(
\int_{\R^N\setminus B_R(x_o)}
\frac{w_\pm^{p-1}(y,t)}{|x-y|^{N+sp}}\,\dy
\right)\dx\,\dt
\notag\\
&\quad\mp
\int_{t_o-S}^{\tau}\int_{B_R(x_o)}
k'(t)\,\xi^p(x)
\left(u^q(x,t)-k^q(t)\right)_\pm\,\dx\,\dt ,
\end{align}
where $c_p>0$ depends only on $p$ and $\gamma>0$ depends only
on the data.
\end{lemma}
\noindent

\section{Positivity estimate} \label{sec5}
\noindent
In this section, we prove a positivity estimate for the solution in terms of moving levels. To this end, we use the following refined version of the Caccioppoli inequality \eqref{cacc2} from Lemma~\ref{cacc}.
\begin{lemma}[Refined energy inequality]\label{caccC}
Let $p>1$, $q>0$, $1<r<R \leq 2$ and $S>0$. Let $u$ be a weak super-solution of \eqref{eq:PDE} in a cylinder containing $\overline B_4 \times [0,S]$, with $u\geq0$ in $B_4\times(0,S]$, with kernel $K=K(x,y,t)$ measurable in the time variable and satisfying \eqref{k}. Let $k \in W^{1,1}(0,S)$ be a function such that $0\leq k(\cdot) \leq \textnormal{\texttt{C}}$, for some $\textnormal{\texttt{C}}>0$. Let $w_\pm(x,t):=(u(x,t)-k(t))_\pm$ and consider a cut-off function $\xi \in C^1_0(B_R; [0,1])$ such that $\xi \equiv 1$ in $B_r$. Define
$$ \textnormal{\texttt{d}} := \textnormal{dist }(\textnormal{supp } \xi, \partial B_R)$$
\begin{equation} \label{PN}
    \mathcal{P}(t):= \int_{\R^N} \frac{u_+^{p-1}(y,t)}{(1+|y|)^{N+sp}}\,\dy, \qquad  \mathcal{N}(t):= \int_{\R^N \setminus B_4} \frac{u_-^{p-1}(y,t)}{|y|^{N+sp}}\,\dy.
\end{equation}
Then, there exist $c_1>0$, $c_2 >1$ depending only on the data such that, for every $\tau \in (0,S]$, it holds
\begin{align} \label{rCC}
   & \int_{B_R} \xi^p \mathfrak{g}_{-}(u(x,\tau),k(\tau))\dx  + \frac{c_p}{\textnormal{\texttt{L}}} \int_0^{\tau} \iint_{B_r \times B_r} \frac{|w_{-}(x,t)-w_{-}(y,t)|^p}{|x-y|^{N+sp}} \dx \dy \dt \notag \\ & \qquad \quad + c_1\int_{0}^{\tau} \int_{B_R} \xi^p w_{-}(x,t) \mathcal{P}(t) \dx \dt \notag \\ &   \qquad \leq \int_{B_R} \xi^p \mathfrak{g}_{-}(u(x,0),k(0))\dx + c_2 \textnormal{\texttt{C}}^p(1+ || \nabla  \xi||_{\infty}^p + \textnormal{\texttt{d}}^{-sp}) |\{ u<k\}\cap B_R \times (0,\tau)| \notag   \\ & \qquad \quad + \int_0^\tau\int_{B_R}\xi^p(x)\left[k'(t)(k^q(t)-u^q(x,t))_+
   +c_2\mathcal N(t)w_-(x,t)\right]\dx\,\dt .
\end{align}
\end{lemma}
\begin{proof}
We start by applying \eqref{cacc2} for supersolutions on $B_R \times (0,\tau)$. Then, for every $\tau \in (0,S]$:
\begin{align} \label{rC}
   & \int_{B_R} \xi^p \mathfrak{g}_{-}(u(x,\tau),k(\tau))\dx - \int_{B_R} \xi^p \mathfrak{g}_{-}(u(x,0),k(0))\dx \notag \\ &  \qquad \quad + \frac{c_p}{\texttt{L}} \int_0^{\tau} \iint_{B_r \times B_r} \frac{|w_{-}(x,t)-w_{-}(y,t)|^p}{|x-y|^{N+sp}} \dx \dy \dt \notag \\ & \qquad \quad + \frac{1}{\texttt{L}} \int_{0}^{\tau} \int_{B_R} \xi^p w_{-}(x,t) \left ( \int_{\R^N} \frac{w_{+}^{p-1}(y,t)}{|x-y|^{N+sp}} \dy \right ) \dx \, \dt \notag \\ &   \qquad \leq \gamma \int_0^{\tau} \iint_{B_R \times B_R} \max \{ w^p_{-}(x,t), w^{p}_{-}(y,t)\} \frac{|\xi(x)-\xi(y)|^p}{|x-y|^{N+sp}} \dx \dy \dt  \notag \\ & \qquad \quad + 4\texttt{L} \int_{0}^{\tau} \int_{B_R} \xi^p w_{-}(x,t) \left ( \int_{\R^N \setminus B_R} \frac{w_{-}^{p-1}(y,t)}{|x-y|^{N+sp}}\,\dy\right ) \dx \dt  \notag  \\ & \qquad \quad +\int_0^\tau k'(t)\int_{B_R}\xi^p(x)(k^q(t)-u^q(x,t))_+\,\dx\,\dt .
\end{align}
Let us estimate the last term on the left-hand side of~\eqref{rC} from below. Note that, by the nonnegativity of $u$ and the bound on $k(\cdot)$, we have
\begin{equation} \label{b1}
   0 \leq w_{-} \leq \texttt{C}, \qquad \text{on } B_R. 
\end{equation}
Now, for $x \in B_2$, it holds 
\begin{equation} \label{b2}
    |x-y| \leq 2(1+|y|).
\end{equation}
Using, $(a+b)^{p-1} \leq 2^{(p-2)_{+}} (a^{p-1}+ b^{p-1})$, $a,b \geq 0$, it follows
\begin{equation} \label{b3}
    w_{+}^{p-1}(y,t) \geq 2^{-(p-2)_{+}} u_{+}^{p-1} (y,t) - \texttt{C}^{p-1}.
\end{equation}
Then, by \eqref{b1}-\eqref{b3}, it holds
\begin{align} \label{b33}
   & \frac{1}{\texttt{L}}\int_{0}^{\tau} \int_{B_R} \xi^p w_{-}(x,t) \left ( \int_{\R^N} \frac{w_{+}^{p-1}(y,t)}{|x-y|^{N+sp}} \dy \right ) \dx \, \dt  \notag \\ & \qquad \quad \geq c_1 \int_{0}^{\tau} \int_{B_R} \xi^p w_{-}(x,t) \mathcal{P}(t)    \dx \dt - \gamma \texttt{C}^{p-1} \int_0^{\tau} \int_{B_R} \xi^p w_{-}(x,t) \dx \dt \notag \\ & \qquad \quad \geq c_1 \int_{0}^{\tau} \int_{B_R} \xi^p w_{-}(x,t)  \mathcal{P}(t) \dx \dt -  \gamma \texttt{C}^{p}|\{ u < k\}\cap B_R \times (0,\tau)|,
\end{align}
with $c_1, \gamma$ depending only  on the data. Denote with $I_1, I_2$ the first two terms on the right-hand side of \eqref{rC}. As regards $I_1$, using \eqref{b1}, we observe that
\begin{align} \label{b4}
    I_1 \leq   \gamma \texttt{C}^p|| \nabla \xi ||_{\infty}^p |\{ u<k\}\cap B_R \times (0,\tau)|,
\end{align}
 with $\gamma$ depending only on the data. For $I_2$,
\begin{align} \label{b5}
I_2 &\leq \gamma \int_0^\tau\!\int_{B_R} \xi^p w_{-}(x,t)
\left[
\texttt{C}^{p-1}\int_{|z|\geq \texttt{d}}\frac{\dz}{|z|^{N+sp}}
+\int_{\mathbb R^N\setminus B_4}
\frac{\texttt{C}^{p-1}+u_-^{p-1}(y,t)}{|y|^{N+sp}}\,\dy
\right]\dx\,\dt  \notag\\
&\leq \gamma(1+\texttt{d}^{-sp})\texttt{C}^p |\{ u<k\}\cap B_R \times (0,\tau)|
+\gamma\int_0^\tau \int_{B_R} \xi^p w_{-}(x,t)
\mathcal{N}(t) \dx\dt 
\end{align}
using  $|x-y| \geq |y|/2$ in $\R^N \setminus B_4$. Therefore, inserting \eqref{b33}-\eqref{b5} into \eqref{rC}, the conclusion \eqref{rCC} follows.
\end{proof}
\noindent
The next lemma gives a positivity estimate in terms of moving levels. Its proof combines the preceding Caccioppoli estimate, applied with moving levels determined by the quantities $\mathcal{P}$ and $\mathcal{N}$ in \eqref{PN}, with a standard De Giorgi iteration. 
\begin{lemma} \label{lem3.2}
Let $\textnormal{\texttt{c}} >0$ be arbitrary. There exists $\delta \in (0,1)$ depending only on the data and $\textnormal{\texttt{c}}$ such that, for $0<S \leq \delta$, a weak super-solution $u$ as in Lemma \ref{caccC} satisfies
\begin{equation} \label{fI}
   u(x,t)\geq\left(h(t)-\ell(t)-\frac{1}{\textnormal{\texttt{c}}}\right)_+^{\frac{1}{q}}
\qquad\text{almost everywhere in } B_1 \times (0,S),
\end{equation}
where
    \begin{equation} \label{hl}
        h(t) = \min \left \{ 2^q, \min \{1,q\} \frac{c_1}{2} \int_0^t \mathcal{P}(\tau) \d \tau  \right \}, \qquad \ell (t) = \max \{1,q\}c_2 \int_{0}^t \mathcal{N}(\tau) \d\tau,
    \end{equation}
  and $\mathcal{P}, \mathcal{N}$ are the quantities defined in \eqref{PN}, $c_1,c_2$ are the constants appearing in Lemma~\ref{caccC}.
\end{lemma}
\begin{proof}
For $n \in \mathbb{N}_0$, we set
    \begin{align} \label{3312}
      &  r_n = 1+2^{-n}, \quad \tilde{r}_n = \frac{r_n + r_{n+1}}{2}, \quad \hat{r}_n = \frac{3 r_n + r_{n+1}}{4}, \quad  B_n = B_{r_n}, \quad \tilde{B}_n = B_{\tilde{r}_n}, \quad \hat{B}_n = B_{\hat{r}_n}, \notag  
     \\ & Q_n := B_n \times (0,S), \quad \tilde{Q}_n := \tilde{B}_n \times (0,S), \quad \hat{Q}_n := \hat{B}_n \times (0,S), \notag \\ & k_n(t) = \left ( h(t)  - \ell (t) - \frac{1}{\texttt{c}}\left(1-\frac{1}{2^n} \right ) \right )_+^{\frac{1}{q}}, \quad w_n = (u-k_n)_{-}, \quad A_n = \{ u < k_n\} \cap Q_n.
    \end{align}
Since \(u,\ell\geq0\), we have
\begin{align}\label{3.12}
&0\leq k_n(t)\leq 2,\qquad
0\leq w_n(x,t)\leq 2
\qquad\text{in }B_4\times(0,S).
\end{align}
Moreover, on \(A_{n+1}\) we have $k_n^q-k_{n+1}^q= \frac{1}{2^{n+1} \texttt{c}}$ and hence
\begin{equation} \label{a2}
    w_n(x,t)\geq k_n(t)-k_{n+1}(t)\geq
\begin{cases}
\displaystyle
\frac{1}{\texttt{c}\,q\,2^{q-1}} \frac{1}{2^{n+1}},
& q\geq1,\\[6pt]
\displaystyle
\left (\frac{1}{\texttt{c}2^{n+1}}\right )^{\frac{1}{q}},
& 0<q<1.
\end{cases}
\end{equation}
 Note that \(
B_{r_{n+1}}\subset \tilde B_n
\subset \hat B_n\subset B_n.
\) We choose $\xi_n \equiv 1$ on $\tilde{B}_n$ and $\text{supp } \xi_n \Subset \hat{B}_{n}$ with $|\nabla \xi_n| \leq 2^{n+4}$. Then $\texttt{d}\geq r_n-\hat r_n=2^{-n-3}$.
 We apply \eqref{rCC} with $(r,R,k,\texttt{C},w_-,\xi)=(\tilde r_n,r_n,k_n,2,w_n,\xi_n)$ and obtain
\begin{align*}
   & \int_{\tilde{B}_n} \mathfrak{g}_{-}(u(x,\tau),k_n(\tau))\dx  + \frac{c_p}{\texttt{L}} \int_0^{\tau} \iint_{\tilde{B}_n \times \tilde{B}_n} \frac{|w_{n}(x,t)-w_{n}(y,t)|^p}{|x-y|^{N+sp}} \dx \dy \dt \notag \\ & \qquad \quad + c_1\int_{0}^{\tau} \int_{B_n} \xi_n^p w_{n}(x,t) \mathcal{P}(t) \dx \dt \notag \\ &   \qquad \leq \int_{{B}_n}  \mathfrak{g}_{-}(u(x,0),k_n(0))\dx + c_2 2^p (1+ ||\nabla \xi_n ||_{\infty}^p + \texttt{d}^{-sp}) |A_n| \notag   \\ & \qquad \quad + \int_0^\tau\int_{B_n}\xi_n^p(x)\big[k_n'(t)(k_n^q(t)-u^q(x,t))_++c_2\mathcal N(t)w_n(x,t)\big]\,\dx\,\dt = I^n_1 + I^n_2 + I^n_3.
\end{align*}
We observe that, when $q>1$, for every $\varepsilon>0$ we can consider $k=k_{n,\varepsilon}:=(k_n-\varepsilon)_+$ that is absolutely continuous and $k_{n,\varepsilon} \to k_n$ as $\varepsilon \to 0$. Let us estimate the terms $I^n_1,I_2^n,I_3^n$ above. We observe that $ I_1^n =0$ since $k_n(0)=0$ and $u$ is nonnegative in $B_4$. As regards $I_2^n$,
\begin{align*}
    I_2^n \leq \gamma 2^p (1+ 2^{pn} + \texttt{d}^{-sp}) |A_n| \leq \gamma  \texttt{b}^n |A_n|.
\end{align*}
For $I_3^n$, we note that $qk_n^{q-1}k_n'=h'-\ell'$ almost everywhere on $\{w_n>0\}$ and
\[
\min \left \{1,\frac{1}{q} \right \}w_n
\leq \frac{k_n^q-u^q}{qk_n^{q-1}}
\leq \max \left \{1,\frac{1}{q} \right \}w_n.
\]
Since $h',\ell'\geq0$, by \eqref{hl} we have
\begin{equation}
 k_n'(k_n^q-u^q)_+ + c_2\mathcal N(t)w_n \leq \begin{cases}
\displaystyle \left(c_2\mathcal N(t)+h'(t)-\frac{\ell'(t)}q\right)w_n \leq \frac{c_1}{2}\mathcal P(t)w_n, & q\geq1,\\[5pt]
\displaystyle \left(c_2\mathcal N(t)+\frac{h'(t)}q-\ell'(t)\right)w_n \leq \frac{c_1}{2}\mathcal P(t)w_n, & 0<q<1.
\end{cases}
\end{equation}
Thus $I_3^n$ can be absorbed in the last term on the left-hand side. Discarding the remaining positive term, taking the supremum on $\tau \in (0,S]$ and using the bounds on $I_1^n, I_2^n$, we get
\begin{align} \label{rCCC}
   \esssup_{\tau \in (0,S]}\int_{\tilde{B}_n} \mathfrak{g}_{-}(u(x,\tau),k_n(\tau))\dx  +  \int_0^{S} \iint_{\tilde{B}_n \times \tilde{B}_n} \frac{|w_{n}(x,t)-w_{n}(y,t)|^p}{|x-y|^{N+sp}} \dx \dy \dt  \leq \gamma \texttt{b}^n |A_n|.
\end{align}

Consider now a spatial cut-off function \(0\leq\phi\leq1\)
with support compactly contained in \(B_{(1-2^{-n-6})\tilde r_n}\), such that \(\phi=1\) on \(B_{r_{n+1}}\)
and 
\(
|D\phi|\leq 2^{n+5}.
\)
Applying Lemma~\ref{em} to \(\phi w_n\), with \(d=2^{-n-6}\) and
\(\mm\) as in \eqref{m}, we obtain
\begin{align}\label{4.9}
&\iint_{\tilde Q_n}
\left[\phi(x,t)w_n(x,t)\right]^{\mm}\,\dx\,\dt
\notag\\
&\quad\leq \gamma\biggl[
\tilde r_n^{sp}
\int_0^S
\iint_{\tilde B_n\times\tilde B_n}
\frac{
\left|\phi(x,t)w_n(x,t)-\phi(y,t)w_n(y,t)\right|^p
}{|x-y|^{N+sp}}\,\dx\,\dy\,\dt
\notag\\
&\qquad\quad
+2^{n(N+sp)}
\iint_{\tilde Q_n}
\left[\phi(x,t)w_n(x,t)\right]^p\,\dx\,\dt
\biggr]
\biggl[
\esssup_{0<t<S}
\fint_{\tilde B_n}
\left[\phi(x,t)w_n(x,t)\right]^{q+1}\,\dx
\biggr]^{\frac{sp}{N}}
\notag\\
&\quad=:\gamma(R_1+R_2)R_3^{\frac{sp}{N}},
\end{align}
where we used \(\tilde r_n\approx 1\). We now estimate the three terms on the right-hand side. Using
\[
\begin{aligned}
&\left|
\phi(x,t)w_n(x,t)-\phi(y,t)w_n(y,t)
\right|^p\\
&\qquad\leq
2^{p-1}|w_n(x,t)-w_n(y,t)|^p\phi^p(x,t)
+2^{p-1}w_n^p(y,t)|\phi(x,t)-\phi(y,t)|^p,
\end{aligned}
\]
together with \eqref{rCCC} and
\eqref{3.12}, we get
\begin{align} \label{R1}
R_1+R_2
&\leq
\gamma\tilde r_n^{sp}
\int_0^S
\iint_{\tilde B_n\times\tilde B_n}
\frac{|w_n(x,t)-w_n(y,t)|^p}{|x-y|^{N+sp}}
\,\dx\,\dy\,\dt
+\gamma\,2^{np}
\iint_{\tilde Q_n}w_n^p(x,t)\,\dx\,\dt
\notag\\
&\qquad
+2^{n(N+sp)}
\iint_{\tilde Q_n}w_n^p(x,t)\,\dx\,\dt
\leq
\gamma\,\texttt{b}^{n}|A_n|.
\end{align}
Set
\[
\beta:=\min\left\{1,\frac{q+1}{2}\right\}.
\]
As regards \(R_3\), we observe that by Lemma \ref{2.3} and \eqref{3.12},
\begin{equation} \label{eqq}
    \begin{cases}
       \frac{1}{q+1} w_n^{q+1}\leq\mathfrak{g}_-(u,k_n) \leq \frac{k_n^{q+1}}{q+1}, \qquad & q\geq 1  \\ \frac{q}{2^{2-q}}w_n^2\leq \mathfrak{g}_-(u,k_n) \leq \frac{k_n^{q+1}}{q+1}, \qquad & 0 < q <1.
    \end{cases}
\end{equation}
Moreover, when $0<q<1$, by H\"older's inequality with exponents \( \left (\frac{2}{q+1}, \frac{2}{1-q} \right )\), we obtain
\begin{equation} \label{snc}
    \fint_{\tilde B_n}(\phi w_n)^{q+1}\,\dx
\leq
\left(\fint_{\tilde B_n}w_n^2\,\dx\right)^{\frac{q+1}{2}},
\end{equation}
since \(0\leq\phi\leq1\). Hence, \eqref{rCCC}, \eqref{eqq}, \eqref{3.12} and \eqref{snc} give
\begin{align}\label{R3}
R_3
&\leq \gamma
\left(
\texttt{b}^{n}|A_n|
\right)^\beta .
\end{align}
Combining \eqref{R1} and \eqref{R3}
with \eqref{4.9}, we obtain
\begin{align}\label{phi}
\iint_{\tilde Q_n}
[\phi(x,t)w_n(x,t)]^\mm\,\dx\,\dt
&\leq
\gamma\texttt{b}^{n}|A_n|
\left(\texttt{b}^{n}|A_n|
\right)^{\frac{sp\beta}{N}}.
\end{align}
Since \(\phi=1\) in \(Q_{n+1}\) and \(w_n=0\) on
\(\tilde Q_n\setminus A_n\), we have
\begin{equation}\label{a1}
\iint_{A_{n+1}}w_n\,\dx\,\dt
\leq
\iint_{\tilde Q_n\cap A_n}\phi w_n\,\dx\,\dt .
\end{equation}
Therefore, by H\"older's inequality with exponents
\(\left(\mm,\frac{\mm}{\mm-1}\right)\), \eqref{a1}, \eqref{a2},
and \eqref{phi}, we get
\begin{align}\label{pe}
\frac{\gamma}{2^{\frac{n}{\min\{1,q\}}}}|A_{n+1}|
&\leq \iint_{A_{n+1}}w_n\,\dx\,\dt
\leq
\left(
\iint_{\tilde Q_n}[\phi(x,t)w_n(x,t)]^\mm\,\dx\,\dt
\right)^{\frac1\mm}
|A_n|^{\frac{\mm-1}{\mm}}
\notag\\
&\leq
\gamma\,\texttt{b}^{n}
|A_n|^{1+\frac{sp\beta}{N\mm}} .
\end{align}
Note that  $\gamma, \texttt{b}$ depend only on the data and on $\texttt{c}$. Define now
$$ Y_n := \frac{|A_n|}{|Q_n|}.$$
Therefore, by \eqref{pe},
it follows that
$$ |A_{n+1}| \leq \gamma\,\texttt{b}^n
|A_n|^{1+\frac{sp\beta}{N\mm}}.$$
Now, we divide by $|Q_{n+1}|$ and use $r_n,r_{n+1}\approx 1$ to get
\begin{align}\label{eq:dg_normalized}
Y_{n+1}
&\leq \gamma\,\texttt{b}^{n}
\frac{|Q_n|}{|Q_{n+1}|}
|Q_n|^{\frac{sp\beta}{N\mm}}
Y_n^{1+\frac{sp\beta}{N\mm}}
\notag\\
&\leq \gamma\,\texttt{b}^{n}
S^{\frac{sp\beta}{N\mm}}
Y_n^{1+\frac{sp\beta}{N\mm}}.
\end{align}
Set \(\alpha:=\frac{sp\beta}{N\mm}\). By the fast geometric convergence
Lemma \ref{fg}, \(Y_n\to0\) as \(n\to\infty\), provided that
\[
Y_0 \leq \,\frac{\nu}{S},
\qquad
\nu:=\gamma^{-1/\alpha}\texttt{b}^{-1/\alpha^2},
\]
which is satisfied taking $S\leq\nu$, since $Y_0\leq1$.
Since $|A_n|\leq Y_n|Q_2|$ for $n\geq2$, $r_n\to1$, and
$k_n(t)\to\left(h(t)-\ell(t)-\frac{1}{\texttt{c}}\right)_+^{\frac1q}$,
we get
$$
\left|\left\{(x,t)\in B_1\times(0,S):
u(x,t)<\left(h(t)-\ell(t)-\frac{1}{\texttt{c}}\right)_+^{\frac1q}
\right\}\right|=0.
$$
Thus,
$$
u(x,t)\geq\left(h(t)-\ell(t)-\frac{1}{\texttt{c}}\right)_+^{\frac{1}{q}}
\qquad\text{almost everywhere in }B_1\times(0,S).
$$
Taking $\delta = \min\{1/2,\nu\}$, the proof is complete.

\end{proof}

\section{Harnack inequalities} \label{sec3}
\noindent
In this section, we prove both the results stated in Theorems~\ref{timegapno}  and~\ref{integralrep}. We start with the Harnack inequality in the fast diffusion regime. The proof is divided into two cases, according to a dichotomy based on the size of the tail. When the tail is small, we follow a local approach. More precisely, we first derive a lower bound for the spatial $L^q$-integral of the solution, which is then used to construct a comparison scheme leading to a positivity estimate.

We then turn to the case of a large tail. The idea is that, when the tail is sufficiently large, the function $h(\cdot)$ in \eqref{hl} becomes large enough to yield a lower bound for the solution through \eqref{fI}. Thus, both cases provide the forward lower bound needed for \eqref{supinf}. The backward upper bound then follows from a classical contradiction argument, see for instance \cite[Section~2]{DGV09}.
\begin{proof}[Proof of Theorem~\ref{timegapno} for $q>p-1>0$
]
We consider the normalized solution
\begin{equation} \label{v}
    v(x,t):=\frac{u(x_o+\varrho x,t_o+u_o^{q+1-p}\varrho^{sp}t)}{u_o},
\qquad u_o:=u(x_o,t_o),
\end{equation}
which solves \eqref{eq:PDE}  in $B_{16} \times (-\mathfrak{a}, \mathfrak{a})$, with the kernel
\[
\tilde K(x,y):=\varrho^{N+sp}K(x_o+\varrho x,x_o+\varrho y),
\]
satisfying \eqref{k} with the same constants. Note that $v(0,0)=1$ and $v\geq0$ globally. Set $Q_{1/16,  \sigma}:= Q_{1/16,  \sigma}(0,0)$. We distinguish two cases according to the size of the  tail:
\begin{align}\label{tail1.1}
\operatorname{Tail}_{p-1}(v;Q_{1/16,  \sigma}  )
&\leq \eps  \sigma^{\frac{N}{sp}},
\end{align}
and
\begin{align}\label{tail1.2}
\operatorname{Tail}_{p-1}(v;Q_{1/16,  \sigma}  )
&> \eps  \sigma^{\frac{N}{sp}},
\end{align}
where $\eps>0$ and $  \sigma>0$ will be chosen later.
\\ \\ \noindent
\textit{Case 1: Assume \eqref{tail1.1}.}\\ \\ \noindent
Let us divide this case into two steps. 
\\ \\ \noindent
\textit{Step 1: A lower bound for the $L^q$-integral of $v$.}\\ \\ \noindent 
We apply the integral Harnack inequality of \cite[Theorem~1.8]{SZ26}. Using $v(0,0)=1$, the nonnegativity of $v$ and \eqref{tail1.1}, we get
\begin{align}\label{hr2}
1\leq{}&
\frac{\gamma}{  \sigma^{\frac{N}{\Lambda_q}}}
\left(\inf_{-  \sigma\leq t\leq0}
\int_{B_{1/4}}v^q(x,t)\,\dx\right)^{\frac{sp}{\Lambda_q}}
+\gamma  \sigma^{\frac1{q+1-p}}
+\gamma \varepsilon^{\frac{sp}{\Lambda_q}},
\end{align}
where $\Lambda_q$ is defined according to \eqref{LAMBDA}. Now, choosing
$\eps,   \sigma\in(0,1)$ such that 
\begin{equation} \label{tildedelta}
    \eps\leq\left (  \frac{1}{4\gamma} \right )^{\frac{\Lambda_q}{sp}}, \qquad   \sigma\leq \left (\frac{1}{4\gamma} \right )^{q+1-p},
\end{equation}
inequality \eqref{hr2} implies
\[
\frac12\leq\frac{\gamma}{  \sigma^{\frac{N}{\Lambda_q}}}
\left(\inf_{-  \sigma\leq t\leq0}
\int_{B_{1/4}}v^q(x,t)\,\dx\right)^{\frac{sp}{\Lambda_q}}.
\]
Therefore, for a constant $\eta_0>0$ depending only on the data,
\begin{align}\label{hr3}
\int_{B_{1/4}}v^q(x,t)\,\dx
\geq \eta_0  \sigma^{\frac{N}{sp}},  
\qquad \textnormal{for all } -  \sigma\leq t\leq0.
\end{align}

\medskip
\noindent
\textit{Step 2: A comparison argument.}\\ \\ \noindent
We introduce $\vartheta:=-  \sigma$ and consider a nonnegative weak solution $w$ of the Cauchy-Dirichlet problem
\begin{equation}\label{DP}
\begin{cases}
\partial_t w^q-\mathcal{L}^s_{\tilde K,p}w=0
&\text{in }B_5\times(\vartheta,  \sigma],\\
w=0
&\text{in }(\R^N\setminus B_5)\times(\vartheta,  \sigma],\\
w^q(\cdot,\vartheta)=\displaystyle
\frac{\eta  \sigma^{\frac{N}{sp}}}{\int_{B_{1/4}}\psi^q(y)\,\dy}
\psi^q
&\text{in }B_5,
\end{cases}
\end{equation}
where $0<\eta<\eta_0$, with $\eta_0$ as in \eqref{hr3}, and
$\psi\in C_c^\infty(B_{1/4})$ is chosen so that
\begin{equation*}
    0\leq\psi\leq v(\cdot,\vartheta),\qquad
\int_{B_{1/4}}\psi^q(y)\,\dy>\eta  \sigma^{\frac{N}{sp}}.
\end{equation*}
Note that the existence of $\psi$ is ensured by the continuity of $v$.  The prescribed initial value for $w$ is a multiple of $\psi$, so existence
follows from the argument of \cite{KMN} that can be easily generalized to a kernel $\tilde K$.
By $(\ref{DP})_3$ and the choice of $\psi$, it holds
\begin{align} \label{4.73}
w(\cdot,\vartheta)
=\left(\frac{\eta  \sigma^{\frac{N}{sp}}}{\int_{B_{1/4}}\psi^q(y)\,\dy}\right)^{\frac1q}
\psi
\leq v(\cdot,\vartheta), \qquad
\int_{B_5}w^q(x,\vartheta)\,\dx
=\eta  \sigma^{\frac{N}{sp}}.
\end{align}
Moreover, since $v$ is globally nonnegative,
\[
w=0\leq v
\qquad\text{a.e. in }(\R^N\setminus B_5)\times(\vartheta,  \sigma].
\]
In particular,
\[
(w-v)_+\in L^p\bigl(\vartheta,  \sigma;W^{s,p}_0(B_5)\bigr).
\]
Since the kernel is independent of time and $w$ is time-independent in $\R^N\setminus B_5$, the comparison principle in \cite{Str26}, applied with $f=0$ after a translation of the time interval, yields
\begin{align}\label{5.27}
w\leq v
\qquad\text{a.e. in }B_5\times(\vartheta,  \sigma].
\end{align}
We now apply the forward $L^q$-$L^q$ estimate in \cite[Remark~1.7]{SZ26}, using $(\ref{DP})_3$. In this estimate no positive tail appears, while the negative tail is zero. Then,
\begin{align}\label{5.28}
\eta  \sigma^{\frac{N}{sp}}
&\leq\gamma\int_{B_{1/2}}w^q(x,t)\,\dx
+\gamma(t-\vartheta)^{\frac{q}{q+1-p}}\leq\gamma\int_{B_{1/2}}w^q(x,t)\,\dx
+\gamma(2  \sigma)^{\frac{q}{q+1-p}},
\qquad \textnormal{for all }\vartheta\leq t\leq  \sigma.
\end{align}
Restricting the value of $  \sigma$ in \eqref{tildedelta} to
\begin{equation} \label{deltadelta}
     \sigma \leq \min \left \{  \left (\frac{1}{4\gamma} \right )^{q+1-p}, \left ( \frac{\eta}{2\gamma} \right )^{\frac{sp(q+1-p)}{\Lambda_q}} \frac{1}{2^{\frac{spq}{\Lambda_q}}}\right \},
\end{equation}
from \eqref{5.28}, we get
\begin{align}\label{5.28a}
\int_{B_{1/2}}w^q(x,t)\,\dx
\geq\frac{\eta}{2\gamma}  \sigma^{\frac{N}{sp}},
\qquad \textnormal{for all }\vartheta\leq t\leq  \sigma.
\end{align}
Now, for every $\tilde{\epsilon}>0$ we consider
\[
\psi_{\tilde{\epsilon}}(\texttt{u})=
\begin{cases}
0, & \texttt{u}\leq 0,\\[2mm]
\dfrac{\texttt{u}}{{\tilde{\epsilon}}}, & 0<\texttt{u}<{\tilde{\epsilon}},\\[2mm]
1, & \texttt{u}\geq{\tilde{\epsilon}}.
\end{cases}
\]
Notice that $0\leq\psi_{\tilde{\epsilon}}\leq1$,
$\psi_{\tilde{\epsilon}}$ is nondecreasing, Lipschitz, and
$\psi_{\tilde{\epsilon}}(0)=0$. By $(\ref{DP})_2$,
$\psi_{\tilde{\epsilon}}(w)$ vanishes outside $B_5$. After the classical regularization in time of $w$ (see \cite[Section 5, formula (5.1)]{CDS}), we choose $\phi = \psi_{\tilde{\epsilon}}(w)$ in 
\eqref{eq} and obtain
\begin{align*}
0=&\int_\vartheta^t\int_{B_5}
\partial_t(w^q)\,\psi_{\tilde{\epsilon}}(w)\,\dx\,\dt\\
&+\int_\vartheta^t
\iint_{\R^N\times\R^N}
|w(x,t)-w(y,t)|^{p-2}
(w(x,t)-w(y,t))
\bigl(\psi_{\tilde{\epsilon}}(w(x,t))
-\psi_{\tilde{\epsilon}}(w(y,t))\bigr)
\tilde K(x,y)\,\dx\,\dy\,\dt.
\end{align*}
By the monotonicity of $\psi_{\tilde{\epsilon}}$, the second term above is nonnegative, then
\[
\int_\vartheta^t\int_{B_5}
\partial_t(w^q)\,\psi_{\tilde{\epsilon}}(w)
\,\dx\,\dt\leq0.
\]
Now, observe that
\[
\Phi_{\tilde{\epsilon}}(\texttt{u})
:=q\int_0^\texttt{u} \xi^{q-1}\psi_{\tilde{\epsilon}}(\xi)\,\d\xi
=\begin{cases}
\displaystyle\frac{q}{q+1}\frac{\texttt{u}^{q+1}}{{\tilde{\epsilon}}},
& 0\leq \texttt{u}\leq{\tilde{\epsilon}},\\[3mm]
\displaystyle \texttt{u}^q-\frac{{\tilde{\epsilon}}^q}{q+1},
& \texttt{u}>{\tilde{\epsilon}},
\end{cases}
\]
therefore $\Phi_{\tilde{\epsilon}}'(\texttt{u})=q\texttt{u}^{q-1}\psi_{\tilde{\epsilon}}(\texttt{u})$
for $\texttt{u}>0$. By the chain rule, justified by time
regularization,
\[
\partial_t(w^q)\,\psi_{\tilde{\epsilon}}(w)
=\partial_t\Phi_{\tilde{\epsilon}}(w),
\]
we get
\begin{align*}
\int_{B_5}\Phi_{\tilde{\epsilon}}(w(x,t))\,\dx
-\int_{B_5}\Phi_{\tilde{\epsilon}}(w(x,\vartheta))\,\dx
&=\int_\vartheta^t\int_{B_5}
\partial_\tau(w^q)\,\psi_{\tilde{\epsilon}}(w)
\,\dx\,\d\tau\leq0.
\end{align*}
Moreover, for $\texttt{u}>{\tilde{\epsilon}}$,
\[
0\leq\Phi_{\tilde{\epsilon}}(\texttt{u})
\leq q\int_0^\texttt{u}\xi^{q-1}\,\d\xi=\texttt{u}^q, \qquad \Phi_{\tilde{\epsilon}}(\texttt{u})=\texttt{u}^q-\frac{{\tilde{\epsilon}}^q}{q+1}
\longrightarrow \texttt{u}^q
\qquad\text{as }{\tilde{\epsilon}} \to 0,     
\]
and $\Phi_{\tilde{\epsilon}}(0)=0$. Hence, by dominated
convergence, using $(\ref{4.73})$, it follows that
\begin{align}\label{5.29a}
\int_{B_5}w^q(x,t)\,\dx
&\leq\int_{B_5}w^q(x,\vartheta)\,\dx
=\eta  \sigma^{\frac{N}{sp}},
\qquad \text{for all }\vartheta\leq t\leq  \sigma.
\end{align}
Now, applying $(\ref{DP})_2$,   H\"older inequality with exponents $\left (\frac{q}{p-1}, \frac{q}{q-p+1} \right )$ and \eqref{5.29a}, we get
\begin{align}\label{5.29b}
\operatorname{Tail}_{p-1}(w;Q_{1/2,  \frac\sigma2}(0,t))\leq \gamma \int_{t-    \frac\sigma2}^{t}\int_{B_5}w^{p-1}(y,\tau)\,\dy\,\d\tau\leq\gamma  \sigma
\left(\eta  \sigma^{\frac{N}{sp}}\right)^{\frac{p-1}{q}},
\qquad -\frac{  \sigma}{2}\leq t\leq\sigma.
\end{align}
Now, applying again the integral Harnack inequality \cite[Theorem~1.8]{SZ26}, using \eqref{5.29a}-\eqref{5.29b} and $\sigma \leq 1$, we get
\begin{align}\label{5.29}
\esssup_{B_{1/2}}w(\cdot,t)
&\leq\frac{\gamma}{  \sigma^{\frac N{\Lambda_q}}}
\left(\eta  \sigma^{\frac N{sp}}\right)^{\frac{sp}{\Lambda_q}}
+\gamma  \sigma^{\frac1{q+1-p}}+\frac{\gamma}{  \sigma^{\frac N{\Lambda_q}}}
\left[  \sigma
\left(\eta  \sigma^{\frac N{sp}}\right)^{\frac{p-1}{q}}\right]^{\frac{sp}{\Lambda_q}}
\leq\tilde\gamma,
\end{align}
for almost every $t\in[-  \sigma/2,  \sigma]$, where $\tilde\gamma>1$ depends only on the data. For the last term we used
\[
-\frac{N}{\Lambda_q}
+\frac{sp}{\Lambda_q}\left(1+\frac{N(p-1)}{spq}\right)=\frac1q.
\]
Then, by \eqref{5.28a} and \eqref{5.29}, we have
\[
\begin{aligned}
\frac{\eta}{2\gamma}  \sigma^{\frac{N}{sp}}
&\leq\int_{B_{1/2}}w^q(x,t)\,\dx\\
&=\int_{B_{1/2}\cap\{w<\mathfrak{b}\}}w^q(x,t)\,\dx
+\int_{B_{1/2}\cap\{w\geq\mathfrak{b}\}}w^q(x,t)\,\dx\\
&\leq\mathfrak{b}^q|B_{1/2}|
+\tilde\gamma^q|\{w(\cdot,t)\geq\mathfrak{b}\}\cap B_{1/2}|,
\end{aligned}
\]
where $\mathfrak{b}>0$ is such that
\begin{equation} \label{FCC}
    \mathfrak{b}\leq \left (\frac{\eta  \sigma^{\frac{N}{sp}}}{4\gamma|B_{1/2}|} \right )^{\frac{1}{q}}.
\end{equation}
Therefore, we obtain
\begin{align}\label{stec1}
|\{w(\cdot,t)\geq\mathfrak{b}\}\cap B_{1/2}|
\geq\alpha|B_{1/2}|,
\qquad -\frac{  \sigma}{2}\leq t\leq\sigma,
\end{align}
where we choose
\[
\alpha:=\frac{\eta  \sigma^{\frac{N}{sp}}}{4\gamma\tilde\gamma^q|B_{1/2}|}\in(0,1).
\]
We now apply the expansion of positivity in \cite{MY}. Note that the tail condition required therein is satisfied by $w$. Then, by \eqref{stec1} and  \cite[Theorem 3.2 and Remark 3.3]{MY}, there exists $\delta^* < \delta_*, \eta_* \in (0,1)$ such that, for every $t_*\in[-\sigma/2,\sigma/4]$,
\begin{align}\label{laut}
w\geq\eta_*\mathfrak{b}
\qquad\text{in }B_1\times
\left(t_*+\delta^*\mathfrak{b}^{q+1-p}\left( \frac{1}{2}\right)^{sp},\,
t_*+\delta_*\mathfrak{b}^{q+1-p}\left( \frac{1}{2}\right)^{sp}\right],
\end{align}
where we further restricted $\mathfrak{b}$ in \eqref{FCC} to
\[
\mathfrak{b}
\leq \min\left\{
\left(\frac{\eta  \sigma^{\frac{N}{sp}}}
{4\gamma|B_{1/2}|}\right)^{\frac{1}{q}},
\left(\frac{  \sigma 2^{sp}}{16\delta^*}\right)^{\frac{1}{q+1-p}}
\right\}.
\]
Then, by the arbitrariness of $t_* \in [-\sigma/2, \sigma/4]$ and  \eqref{5.27}, we obtain
\begin{align}\label{ehr}
v\geq w\geq\eta_*\mathfrak{b}
\qquad\text{in }B_1\times
\left[-\frac{  \sigma}{8},\frac{  \sigma}{8}\right].
\end{align}
This concludes \textit{Case 1}.

\medskip
\noindent
\textit{Case 2: Assume \eqref{tail1.2}.}\\ \\ \noindent
We observe that, by \eqref{tail1.2} and the definition of $\mathcal P$ in \eqref{PN}, we have
\begin{align}\label{cpt}
\int_{-  \sigma}^0\mathcal P_v(t)\,\dt
&\geq\frac{1}{17^{N+sp}}\operatorname{Tail}_{p-1}(v;Q_{1/16,  \sigma}  )>\frac{\eps  \sigma^{\frac{N}{sp}}}{17^{N+sp}}.
\end{align}
Now, let $0<\lambda\leq1$ be a constant to be chosen later. We introduce
\[
z(x,\tau):=\frac{v(x,-  \sigma+\lambda^{q+1-p}\tau)}{\lambda}
\]
and set
\begin{align}\label{def:tau_o}
\tau_o:=\frac{  \sigma}{\lambda^{q+1-p}}.
\end{align}
Note that $z$ is a globally nonnegative weak solution to \eqref{eq:PDE}, with kernel $\widetilde K$, in a cylinder containing $B_4\times[0,1]$. In Lemma~\ref{lem3.2}, we choose
\begin{equation}\label{cq}
\texttt{c}=c_q:=2\frac{1+\max\{1,q\}}{2^q-1},
\end{equation}
and let $\delta$ denote the corresponding constant of that lemma. We apply the lemma on $B_4\times(0,S)$, with $S\in(0,\delta]$ to be chosen later. It yields
\begin{equation}\label{suV}
z(x,\tau)\geq\left(h(\tau)-\ell(\tau)-\frac1{c_q}\right)_+^{\frac1q}
\end{equation}
for almost every $(x,\tau)\in B_1\times(0,S)$, where $h,\ell$ are defined by \eqref{hl} with $u=z$.
By \eqref{cpt} and the change of variables,
\begin{align}\label{suH}
\int_0^{\tau_o}\mathcal P_z(\tau)\,\d\tau
&=\frac1{\lambda^q}\int_{-  \sigma}^0\mathcal P_v(t)\,\dt>\frac{\eps}{17^{N+sp}}\frac{  \sigma^{\frac N{sp}}}{\lambda^q}.
\end{align}
Consequently, provided that $\tau_o<S$, for all $\tau\in[\tau_o,S]$ it holds
\begin{equation} \label{H2}
    h(\tau)\geq\min\left\{2^q,
\frac{c_1\min\{1,q\}}2\,\frac{\eps}{17^{N+sp}}
\frac{  \sigma^{\frac N{sp}}}{\lambda^q}\right\}.
\end{equation}
On the other hand, since $v\geq0$ globally, we have
\begin{equation}\label{suL}
\ell(\tau)=0,\qquad \tau\in(0,S).
\end{equation}
We choose
\begin{align}\label{c4c5}
\lambda:=\mu  \sigma^{\frac{N}{spq}} \quad \textnormal{with}\quad  \mu\in(0,1) \quad \text{such that} \quad \frac{c_1\min\{1,q\}}{2}\,\frac{1}{17^{N+sp}}\frac{\eps}{\mu^q}\geq2^q.
\end{align}
Moreover, since $\Lambda_q>0$, we can further reduce $  \sigma$ so that
\[
2\tau_o=2\mu^{-(q+1-p)}
  \sigma^{1-\frac{N(q+1-p)}{spq}}
=2\mu^{-(q+1-p)}
  \sigma^{\frac{\Lambda_q}{spq}}
\leq\delta.
\]
We then set $S=2\tau_o$. Now, by \eqref{H2},  \eqref{suL} and \eqref{c4c5}, we get
\begin{equation}\label{hll}
h(\tau)=2^q,\qquad\ell(\tau)=0,
\qquad \tau_o<\tau<2\tau_o.
\end{equation}
Consequently, by \eqref{suV} and \eqref{hll}, it follows that
\[
z(x,\tau)\geq\left(2^q-\frac1{c_q}\right)^{\frac{1}{q}}>0,
\qquad\text{a.e. in }B_1\times(\tau_o,2\tau_o).
\]
Returning to $v$, we obtain
\begin{align}\label{concl1}
v(x,t)\geq
\mu\left(2^q-\frac1{c_q}\right)^{\frac{1}{q}}
  \sigma^{\frac{N}{spq}},
\qquad \text{a.e. in }B_1\times(0,  \sigma).
\end{align}
This finishes \textit{Case 2}.

\medskip
\noindent
\textit{Conclusion: Harnack inequality.}\\ \\ \noindent
We set
\begin{equation} \label{deltaacca}
    \sigma_{\texttt{f}}:=\frac{  \sigma}{8},\qquad
\gamma_\texttt{f}:=\left[\min\left\{
\eta_*\mathfrak{b},\,
\mu\left(2^q-\frac1{c_q}\right)^{\frac{1}{q}}
  \sigma^{\frac{N}{spq}},\,\frac12\right\}\right]^{-1}.
\end{equation}
Combining \eqref{ehr} and \eqref{concl1}, and returning to the original variables, we conclude that
\begin{equation}\label{3.34}
u_o\leq\gamma_\texttt{f}\inf_{B_\varrho(x_o)}u(\cdot,t),
\qquad\text{for all }t\in
\left(t_o,\,t_o+\sigma_{\texttt{f}}u_o^{q+1-p}\varrho^{sp}\right),
\end{equation}
provided
\[
B_{6\varrho}(x_o)\times
\left[t_o-  \sigma u_o^{q+1-p}\varrho^{sp},\,
t_o+  \sigma u_o^{q+1-p}\varrho^{sp}\right]\Subset\Omega_T.
\]
Next, we claim that
\begin{equation} \label{sup}
    \sup_{B_\varrho(x_o)}u(\cdot,t)
    \leq 2\gamma_\texttt{f}u_o, \qquad \text{for all }  t\in
    \left(
        t_o-\sigma_{\texttt{f}}u_o^{q+1-p}\varrho^{sp},
        \,
        t_o
    \right).
\end{equation}
Indeed, suppose that the contrary of
\eqref{sup} holds.
Since $u(x_o,t_o)=u_o$, by continuity, there exists
\begin{equation}
    (x_*,t_*)\in
    B_\varrho(x_o)\times
    \left(
        t_o-\sigma_{\texttt{f}}u_o^{q+1-p}\varrho^{sp},
        \,
        t_o
    \right)
    \label{star}
\end{equation}
such that
\begin{equation}
   u_*:= u(x_*,t_*)=2\gamma_\texttt{f}u_o.
    \label{starV}
\end{equation}
Now, let us observe that, since $x_*\in B_\varrho(x_o)$,
\begin{equation} \label{suB}
    B_{12\varrho}(x_*)
\Subset B_{16\varrho}(x_o).
\end{equation}
Moreover, by \eqref{star} and
\eqref{starV},
\begin{equation}\label{suT}
\left[
    t_*-  \sigma u_{*}^{q+1-p}(2\varrho)^{sp},
    \,
    t_*+  \sigma u_{*}^{q+1-p}(2\varrho)^{sp}
\right]\subset
\left(
    t_o-\mathfrak{a} u_o^{q+1-p}\varrho^{sp},
    \,
    t_o+\mathfrak{a} u_o^{q+1-p}\varrho^{sp}
\right),
\end{equation} 
provided
\begin{equation}
    \label{aio}
    \mathfrak{a}>
    \sigma_{\texttt{f}}+  \sigma(2\gamma_\texttt{f})^{q+1-p}2^{sp}.
\end{equation}
Then, due to \eqref{suB} and \eqref{suT}, we can apply \eqref{3.34} at
$(x_*,t_*)$ and obtain
\begin{equation}
    2\gamma_\texttt{f}u_o
    =
    u_{*}
    \leq
    \gamma_\texttt{f}\inf_{B_{2\varrho}(x_*)}u(\cdot,t), \qquad \text{for all }  t\in
    \left(
        t_*,
        \,
        t_*+\sigma_{\texttt{f}}u_{*}^{q+1-p}
        (2\varrho)^{sp}
    \right).
    \label{suu}
\end{equation}
Now, by \eqref{star} and
\eqref{starV},
\[
    t_o-t_*
    <\sigma_{\texttt{f}}u_o^{q+1-p}\varrho^{sp}
    <
    \sigma_{\texttt{f}}u_{*}^{q+1-p}
    (2\varrho)^{sp}.
\]
Since $x_o\in B_{2\varrho}(x_*)$, it follows that
\[
    (x_o,t_o)\in
    B_{2\varrho}(x_*)\times
    \left(
        t_*,
        \,
        t_*+\sigma_{\texttt{f}}u_{*}^{q+1-p}
        (2\varrho)^{sp}
    \right).
\]
Taking $t=t_o$ in \eqref{suu},
we obtain
\[
    2\gamma_\texttt{f}u_o\leq\gamma_\texttt{f}u_o,
\]
a contradiction. Therefore
\eqref{sup} holds. Combining \eqref{sup} with
\eqref{3.34}, we conclude that
\[
    \sup_{ B_\varrho(x_o)\times
    \left (
        t_o-\sigma_{\texttt{f}}u_o^{q+1-p}\varrho^{sp},
        \,
        t_o \right )}u
    \leq 2\gamma_\texttt{f}^2\inf_{B_\varrho(x_o)\times
    \left (
        t_o,
        \,
        t_o+\sigma_{\texttt{f}}u_o^{q+1-p}\varrho^{sp} \right )}u,
\]
which implies the conclusion \eqref{supinf} in the fast diffusion regime.
\end{proof}
\noindent
We turn to the proof of the slow diffusion case. The proof incorporates also $q=p-1$. Starting from a quantitative local upper estimate we first derive an integral representation estimate. Combining this with the moving-level estimate \eqref{fI} we obtain the lower bound for $u$ used to prove the upper bound in \eqref{supinf}. The conclusion follows again by a contradiction argument.

\begin{proof}[Proof of Theorem \ref{timegapno} for $0<q \leq p-1$.]
We prove the more general formulation under \eqref{negativetaill}.
Set
\begin{equation} \label{b}
    b:= \begin{cases}
        1, \qquad \text{if } q<p-1
        \\ \frac{\delta}{2}, \qquad \text{if } q=p-1
    \end{cases}
\end{equation}
where $\delta$ is as in Lemma~\ref{lem3.2} for the choice of $\texttt{c}$ as in~\eqref{cq}. We use the local upper bound from \cite[Proposition~1.4]{SSY2026_slow}. To extend this estimate to the full range $0<q\leq p-1$, we observe that the assumption \cite[(5.3)]{SSY2026_slow} is not essential. Indeed, the corresponding choice of $k$ becomes ill-defined at $q=p-1$, as its defining expression may either diverge or tend to zero. However, by suitably adapting the cylinders to the intrinsic scaling, the estimate remains valid also in the limiting case $q=p-1$, with constants depending only on the data.

Let $\tilde{\texttt{a}}\geq1$ denote the corresponding constant and set
\begin{equation}\label{slow:theta}
\theta:=b\left(\frac{u_o}{4\tilde{\texttt{a}}}\right)^{q+1-p}\varrho^{sp},\qquad u_o:=u(x_o,t_o).
\end{equation}
Take $\mathfrak{a}$ as
\begin{align}\label{aio2}
\mathfrak{a}>(4\tilde{\texttt{a}})^{p-q-1}
,
\end{align}
we can apply this estimate to the cylinder 
\begin{align*}
Q^-_{\varrho, \theta}(x_o,t_o)
=
B_\varrho(x_o)\times(t_o- \theta,t_o).
\end{align*}
Thus,
\begin{align*}
u_o
&\leq
\operatorname*{ess\,sup}_{
Q^-_{\varrho/2, \theta/2}(x_o,t_o)}u \leq
\tilde{\texttt{a}}
\left(
\int_{t_o- \theta}^{t_o}
\int_{\mathbb R^N}
\frac{u_+^{p-1}(x,t)}
{(\varrho+|x-x_o|)^{N+sp}}
\, \dx\, \dt
\right)^{\frac{1}{q}}
+\frac{1}{4 } u_o.
\end{align*}
From which, it follows that
\begin{align} \label{usu}
    u_{o}^{q} \leq ( 2 \tilde{\texttt{a}} )^{q} \int_{t_o-\theta}^{t_o}
\int_{\mathbb R^N}
\frac{u_+^{p-1}(x,t)}
{(\varrho+|x-x_o|)^{N+sp}}
\, \dx\, \dt.
\end{align}
 Now, we use this upper bound to prove \eqref{supinf}. Let $\lambda \in (0, u_o]$ be a constant that will be determined later and define
\begin{align*}
    \widetilde{v}(z, \tau) := \frac{u(x_o + \varrho z, t_o - \theta + \lambda^{q+1-p} \varrho^{sp} \tau )}{\lambda},
\end{align*}
which is a weak solution to~\eqref{eq:PDE}, where the kernel satisfies~\eqref{k}. Now,  let $\texttt{c}$ be as in \eqref{cq} and $\delta\in(0,1)$ be the corresponding number in Lemma \ref{lem3.2}. Consider $\theta$ as \eqref{slow:theta} and $S=2\tau_o$, where $\tau_o:=\theta/(\lambda^{q+1-p}\varrho^{sp})$. A correct choice of $\lambda$ will indeed show $S \leq \delta$.
Now, the cylinder $B_4\times[0,S]$ corresponds in the original
variables to $B_{4\varrho}(x_o)\times[t_o-\theta,t_o+\theta]$,
which is compactly contained in $\Omega_T$ by \eqref{aio2}.
Hence, the cylinder required by Lemma~\ref{lem3.2} is admissible.
By applying Lemma~\ref{lem3.2}, we deduce for $\widetilde{v}(z, \tau)$ the estimate in~\eqref{fI} for a.e. $(z,\tau)\in B_1\times(0,S)$, where $h,\ell$ are defined by \eqref{hl} for the rescaled solution $\widetilde v$. The estimate in~\eqref{usu} combined with a change of variables yields
\begin{align*}
\int_0^{\tau_o}\mathcal{P}_{\widetilde{v}}(\tau)\,\d\tau
&=\frac{1}{\lambda^q}
\int_{t_o-\theta}^{t_o}
\int_{\mathbb R^N}
\frac{u_+^{p-1}(x,t)}
{(\varrho+|x-x_o|)^{N+sp}}
\, \dx\, \dt \geq \frac{1}{(2\tilde{\texttt{a}})^q} \Big( \frac{u_o}{\lambda} \Big)^{q},
\end{align*}
As a consequence, for all $\tau\in[\tau_o,S]$, observe that $\tau_o <S$, it follows that
\begin{align}\label{suH_two}
h(\tau)\geq
\min\left\{
2^q,\,
\frac{c_1\min\{1,q\}}{2}
\frac{1}
{(2\tilde{\texttt{a}})^q}
\left(\frac{u_o}{\lambda}\right)^q
\right\}.
\end{align}
Besides, by \eqref{negativetaill} and a change of variables,
for any $\tau \in (0,S)$, we deduce
\[
\ell(\tau)\leq \max\{1,q\}c_2\mathfrak{m}
\left(\frac{u_o}{\lambda}\right)^q,\qquad \tau\in(0,S).
\]
Next, let $\lambda:=
\frac{\zeta}{4 \tilde{\texttt{a}}} u_o$ with $\zeta \in (0,1)$  chosen in a way such that
\begin{align}\label{def:k_two}
    \frac{c_1\min\{1,q\}}{2\zeta^q}\geq 2^q,\qquad 2b\zeta^{p-(q+1)}\leq\delta,
\end{align}
note that the second choice is possible by \eqref{b} and implies $S \leq \delta$.
Moreover, take the constant $\mathfrak{m} \in (0,1)$ in \eqref{negativetaill} sufficiently small such that
\begin{equation} \label{eta_two}
    \max\{1,q\}c_2\mathfrak{m}
\left(\frac{4\tilde{\texttt{a}}}{\zeta}\right)^q
\leq \frac{1}{c_q}.
\end{equation}
Combining~\eqref{suH_two},~\eqref{def:k_two} and~\eqref{eta_two} yields
\begin{align*}
    \widetilde{v}(z,\tau)\geq
\left(2^q-\frac{2}{c_q}\right)^{\frac{1}{q}}>0
\end{align*}
for a.e. $(z,\tau)\in B_1\times(\tau_o,2\tau_o)$. For the initial variables this means that
\begin{align*}
    u(x,t)\geq \frac{\zeta}{4\tilde{\texttt{a}}}
\left(2^q-\frac{2}{c_q}\right)^{\frac{1}{q}} u_o, \qquad \text{a.e. in } B_\varrho(x_o)
\times(t_o,t_o+ \theta).
\end{align*}
Now, let us define the two quantities
\begin{align} \label{sigmaesse}
 \gamma_\texttt{s}:=
\frac{4\tilde{\texttt{a}}}{\zeta}
\left(2^q-\frac{2}{c_q}\right)^{-\frac{1}{q}}, \qquad \sigma_{\texttt{s}}:=b\min\left\{
\frac12,\,\frac{\zeta}{(4\tilde{\texttt{a}})^{q+1-p}},
\frac12(4\tilde{\texttt{a}})^{p-q-1}
(2\gamma_\texttt{s})^{q+1-p}2^{sp}
\right\}.
\end{align}
By the continuity we deduce
\begin{align}\label{ineq_u_o_le_inf_u_estimate}
    u_o\leq \gamma_\texttt{s}
\inf_{B_\varrho(x_o)}u(\cdot,t),
\qquad \text{for all }
t\in\left(
t_o,\,
t_o+\sigma_{\texttt{s}} u_o^{q+1-p}\varrho^{sp}
\right)
\end{align}
and claim that
\begin{align}\label{ineq:sup_claim}
    \sup_{B_\varrho(x_o)}u(\cdot,t)
    \leq 2\gamma_\texttt{s}u_o, \qquad \text{for all }  t\in
    \left(
        t_o-\sigma_{\texttt{s}}u_o^{q+1-p}\varrho^{sp},
        \,
        t_o
    \right).
\end{align}
Enlarging $\mathfrak{a}$ in \eqref{aio2} to
\begin{align}\label{aio4}
\mathfrak{a}>\max\left\{
(4\tilde{\texttt{a}})^{p-q-1},
\sigma_{\texttt{s}}+(4\tilde{\texttt{a}})^{p-q-1}
(2\gamma_\texttt{s})^{q+1-p}2^{sp}
\right\},
\end{align}
the claim follows by the same argument as in the proof
of \eqref{sup}, with the negative tail controlled by
\eqref{negativetaill}. Indeed, setting
\[
\mathfrak{m}_*:=\left(\frac32\right)^{N+sp}
\frac{\mathfrak{m}}{(2\gamma_\texttt{s})^q},
\]
we further restrict $\mathfrak{m}$ so that
\[
\max\{1,q\}c_2\mathfrak{m}_*
\left(\frac{4\tilde{\texttt{a}}}{\zeta}\right)^q
\leq \frac{1}{c_q}.
\]
The choice of $\sigma_{\texttt{s}}$ ensures that $t_o$
belongs to the forward positivity interval associated
with the new center. Thus,~\eqref{ineq:sup_claim} combined with~\eqref{ineq_u_o_le_inf_u_estimate} yields 

\begin{align} \label{supinfslow}
    \sup_{ B_\varrho(x_o)\times
    \left (
        t_o-\sigma_{\texttt{s}}u_o^{q+1-p}\varrho^{sp},
        \,
        t_o \right )}u
    \leq 2\gamma_\texttt{s}^2\inf_{B_\varrho(x_o)\times
    \left (
        t_o,
        \,
        t_o+\sigma_{\texttt{s}}u_o^{q+1-p}\varrho^{sp} \right )}u,
\end{align}
which implies \eqref{supinf} in the slow diffusion regime. 
Note that the structural constants $\mathfrak{a}$, $\mathfrak{s}$
and $\mathfrak{c}$ may differ in the two cases $q>p-1$
and $0<q\leq p-1$. The proof of Theorem~\ref{timegapno} is complete.
\end{proof}
\begin{remark}[Integral estimates in the slow-diffusion range]
The estimate~\eqref{usu} holds both for $sp<N$ and for
$N\leq sp$. In the latter case, it suffices to replace
the Sobolev embedding from~\cite[Lemma 2.1]{SSY2026_slow}
with that from~\cite[Lemma 2.2]{SZ26}
in the proof of~\cite[Proposition 1.4]{SSY2026_slow}.
\end{remark}

\noindent
We finally turn to the integral representation estimates. For the lower bound, the idea is that if the positive contribution of $h$ in \eqref{hl} becomes sufficiently large, while the negative tail in $\ell$ remains small, the positivity estimate forces the solution to exceed its value at the reference point, leading to a contradiction. This follows the approach developed by Liao in \cite{NN}. As for the upper bound, in the slow diffusion regime (and in the nonlocal Trudinger equation) it follows again from a local $L^\infty$-estimate, whereas in the fast diffusion regime it is derived from the Harnack inequality established in \eqref{supinf}.
\begin{proof}[Proof of Theorem~\ref{integralrep}]
Let us start proving the lower bound in \eqref{integralbounds}. Let $\delta$ be the constant in Lemma~\ref{lem3.2}
corresponding to the choice $\texttt{c}$ in~\eqref{cq}. Let $v$ be the normalized solution in \eqref{v}. In particular, $v(0,0)=1$ and $v\geq0$
globally. Let $\mathcal P_v, \mathcal{N}_v$ is defined as in \eqref{PN} with $u=v$. We apply
Lemma~\ref{lem3.2} to $v$ after translating the initial time
$-{\frac\delta2}$ to $0$. We have
\begin{equation}\label{rep-hl}
 h(t)=\min\left\{2^q,\min\{1,q\}\frac{c_1}{2}
 \int_{-{\frac\delta2}}^t\mathcal P_v(\tau)\,\d\tau\right\},
 \qquad
 \ell(t)=\max\{1,q\}c_2
 \int_{-{\frac\delta2}}^t\mathcal N_v(\tau)\,\d\tau.
\end{equation}
Since $v\geq0$ globally, $\mathcal N_v=0$ and therefore
$\ell=0$. By Lemma~\ref{lem3.2} and the continuity of $v$ and $h$,
\begin{equation}\label{rep-point}
1=v(0,0)\geq\left(h(0)-\ell(0)-\frac1{c_q}\right)_+^{\frac{1}{q}}.
\end{equation}
Thus, if
\begin{equation*}
   \int_{-{\frac\delta2}}^0\mathcal P_v(\tau)\,\d\tau
\geq\frac{2^{q+1}4}{\min\{1,q\}c_1}, 
\end{equation*}
it follows that $h(0)=2^q$ and so by \eqref{rep-point} and
\eqref{cq} we have
\[
1\geq\left(2^q-\frac1{c_q}\right)^{\frac{1}{q}}
\geq\left(\frac{2^q+1}{2}\right)^{\frac{1}{q}}>1,
\]
a contradiction. Therefore,
$$ \int_{-{\frac\delta2}}^0 \mathcal{P}_v(\tau) \d \tau <  \frac{2^{q+1}4}{\min\{1,q\}c_1},$$
and  by the change of variables 
$$\frac{1}{u_o^q} \int_{t_o - {\frac\delta2}u_o^{q+1-p}\rho^{sp}}^{t_o} \int_{\R^N} \frac{u_{+}^{p-1}(x,t)}{(\rho + |x-x_o|)^{N+sp}} \dx \dt = \int_{-{\frac\delta2}}^0 \mathcal{P}_v(\tau) \d \tau \leq \frac{2^{q+1}4}{\min\{1,q\}c_1},$$
which is the lower bound in \eqref{integralbounds}. 

Now, let us move to the upper bound in \eqref{integralbounds2}. In the regime $q\leq p-1$, we repeat the argument leading to
\eqref{usu} with 
\[
r=\kappa\varrho, \qquad \kappa:=
\left(
\frac{\delta}{2(4\tilde{\texttt{a}})^{p-q-1}}
\right)^{\frac1{sp}}, \qquad \theta_r:=
\left(\frac{u_o}{4\tilde{\texttt{a}}}\right)^{q+1-p}r^{sp}
={\frac\delta2}u_o^{q+1-p}\varrho^{sp}.
\]
Applying the same local upper estimate on the cylinder
\[
Q^-_{r,\theta_r}(x_o,t_o)
=B_r(x_o)\times(t_o-\theta_r,t_o),
\]
we obtain
\begin{equation} \label{rr}
    u_o^q\leq (2\tilde{\texttt{a}})^q
\int_{t_o-\theta_r}^{t_o}
\int_{\mathbb R^N}
\frac{u_+^{p-1}(x,t)}
{(r+|x-x_o|)^{N+sp}}
\,\dx\,\dt.
\end{equation}
Since
\(
r+|x-x_o|\geq
\kappa(\varrho+|x-x_o|),
\)
the upper bound \eqref{integralbounds2} follows,
with a constant depending only on the data.

Let us analyze the fast diffusion case:  $q>p-1$. Let $\mathfrak{s},  \mathfrak{c}, \mathfrak{a}$ be as in the statement of
Theorem~\ref{timegapno}. Set
\[
\mathcal J:=\int_{-\mathfrak{s}}^0\int_{\R^N}
\frac{v^{p-1}(x,t)}{(1+|x|)^{N+sp}}\,\dx\,\dt.
\]Applying \eqref{supinf} to
the normalized solution $v$ and also due to continuity, we deduce
\[
\sup_{B_1\times(-{\mathfrak{s}},0)}v
\leq\mathfrak{c}\inf_{B_1\times(0,{\mathfrak{s}})}v
\leq\mathfrak{c} v(0,0)=\mathfrak{c}.
\]
We observe that
\[
v^q=v^{p-1}v^{q+1-p}\leq \mathfrak{c}^{q+1-p}v^{p-1}
\qquad\text{in }B_1\times(-{\mathfrak{s}},0),
\]
hence,
\begin{align}\label{rep-mass}
\inf_{-{\mathfrak{s}}\leq t\leq0}\int_{B_{1/4}}v^q(x,t)\,\dx
&\leq\frac1{{\mathfrak{s}}}\int_{-{\mathfrak{s}}}^0\int_{B_{1/4}}v^q(x,t)\,\dx\,\dt\leq\frac{\gamma}{{\mathfrak{s}}}\mathcal J,
\end{align}
with $\gamma$ depending on data. Moreover, we observe that
\begin{equation}\label{rep-tail}
\operatorname{Tail}_{p-1}(v;Q_{1/16,\mathfrak{s}}  )
\leq \gamma\mathcal J.
\end{equation}       
We apply the integral Harnack inequality
\cite[Theorem~1.8]{SZ26}, by \eqref{rep-mass} and \eqref{rep-tail}, it follows that
\begin{align*}
1&\leq\frac{\gamma}{\mathfrak{s}^{\frac N{\Lambda_q}}}
\left(\inf_{-\mathfrak{s}\leq t\leq0}\int_{B_{1/4}}v^q(x,t)\,\dx\right)^{\frac{sp}{\Lambda_q}}
+\gamma\mathfrak{s}^{\frac1{q+1-p}}
+\frac{\gamma}{\mathfrak{s}^{\frac N{\Lambda_q}}}
\operatorname{Tail}_{p-1}(v;Q_{1/16,\mathfrak{s}}  )^{\frac{sp}{\Lambda_q}}\\
&\leq\gamma\mathfrak{s}^{-\frac{N+sp}{\Lambda_q}}\mathcal J^{\frac{sp}{\Lambda_q}}
+\gamma\mathfrak{s}^{\frac1{q+1-p}}
+\gamma\mathfrak{s}^{-\frac N{\Lambda_q}}\mathcal J^{\frac{sp}{\Lambda_q}}\\
&\leq\gamma\mathfrak{s}^{-\frac{N+sp}{\Lambda_q}}\mathcal J^{\frac{sp}{\Lambda_q}}
+\gamma\mathfrak{s}^{\frac1{q+1-p}}.
\end{align*}
Take $\mathfrak{s} \leq  \left (\frac{1}{4\gamma} \right )^{q+1-p}$ then, the last term is at most $1/4$.
Thus,
\[
\frac34\leq\gamma\mathfrak{s}^{-\frac{N+sp}{\Lambda_q}}
\mathcal J^{\frac{sp}{\Lambda_q}},
\]
which, by the change of variables, gives
\begin{equation}\label{rep-upper}
\frac{1}{u_o^q}
\int_{t_o-\mathfrak{s} u_o^{q+1-p}\varrho^{sp}}^{t_o}
\int_{\R^N}
\frac{u^{p-1}(x,t)}{(\varrho+|x-x_o|)^{N+sp}}\,\dx\,\dt
= \mathcal J\geq\frac1\gamma\mathfrak{s}^{\frac{N+sp}{sp}}
\geq\frac1{\tilde\gamma},
\end{equation}
with $\tilde\gamma>1$ depending only on the data. Taking $\mathfrak{s}\leq \frac\delta2$ and \(\mathfrak c\) sufficiently large complete the proof.
\end{proof}
\noindent
The following remark explains how to obtain the lower bound when the solution is only locally nonnegative.
\begin{remark}[Local nonnegativity in the integral lower bound]\label{1rem}
We observe that the global nonnegativity assumption in the proof of the
lower bound in \eqref{integralbounds} can be replaced by
local nonnegativity and \eqref{rep-negative-tail-condition}. Moreover, we can consider a kernel $K=K(x,y,t)$, measurable in the time variable.
Indeed, let $\mathfrak{s}$ be as in the statement of Theorem \ref{integralrep}, note that  by the proof above $\mathfrak{s} \leq \delta/2$. Then, with $v$ that is the normalized solution in \eqref{v}, by \eqref{rep-negative-tail-condition}, we get
\[
\int_{-\mathfrak{s}}^0\mathcal N_v(t)\,\dt
=\frac1{u_o^q}
\int_{t_o-\mathfrak{s}u_o^{q+1-p}\varrho^{sp}}^{t_o}
\int_{\R^N\setminus B_{4\varrho}(x_o)}
\frac{u_-^{p-1}(x,t)}{|x-x_o|^{N+sp}}\,\dx\,\dt
\leq\mathfrak{m}.
\]
Let $h, \ell$ be as in \eqref{hl},  after translating
the initial time $-\mathfrak{s}$ to $0$. Choosing $\mathfrak{m} \leq \frac{1}{c_q c_2},$ we get
\[
\ell(0)\leq \max\{1,q\}c_2 \mathfrak{m}
\leq\frac{\max\{1,q\}}{c_q}.
\]
Assuming
\begin{equation*}
   \int_{-\mathfrak{s}}^0\mathcal P_v(t)\,\dt
\geq\frac{2^{q+1}4}{\min\{1,q\}c_1}, 
\end{equation*}
we obtain $h(0)=2^q$. Then, Lemma~\ref{lem3.2} and \eqref{cq}, imply
\[
1=v(0,0)\geq
\left(2^q-\frac{1+\max\{1,q\}}{c_q}\right)^{\frac{1}{q}}
=\left(\frac{2^q+1}{2}\right)^{\frac{1}{q}}>1,
\]
a contradiction. Consequently,
\[
\int_{-\mathfrak{s}}^0\mathcal P_v(t)\,\dt
<\frac{2^{q+1}4}{\min\{1,q\}c_1},
\]
which, after rescaling, gives \eqref{integralbounds},
with the positive part $u_+^{p-1}$ in the integral.
\end{remark}
\noindent
We next show that the Harnack inequality of Theorem~\ref{timegapno}
can be strengthened by including the positive tail on the
left-hand side.
\begin{theorem}[Full Harnack inequality without time gap]\label{fullharnack}
Under the assumptions of Theorem \ref{timegapno}, there exist
constants $\mathfrak{s}\in(0,1)$, $\mathfrak{c}>1$ and
$\mathfrak{a}>1$, depending only on the data, such that if
\[
B_{16\varrho}(x_o)\times
\left[
t_o-\mathfrak{a}u_o^{q+1-p}\varrho^{sp},
t_o+\mathfrak{a}u_o^{q+1-p}\varrho^{sp}
\right]\Subset\Omega_T
\]
then
\begin{equation}\label{fullsupinf}
\sup_{Q^-}u+
\textnormal{Tail}_{p-1}^{\frac{1}{q}}(u;Q^-)
\leq \mathfrak{c}\inf_{Q^+}u,
\end{equation}
where
\[
Q^-:=B_\varrho(x_o)\times
(t_o-\mathfrak{s}u_o^{q+1-p}\varrho^{sp},t_o],
\qquad
Q^+:=B_\varrho(x_o)\times
(t_o,t_o+\mathfrak{s}u_o^{q+1-p}\varrho^{sp}].
\]
\end{theorem}
\noindent
The proof follows straightforwardly by combining \eqref{supinf} and
\eqref{integralbounds}.

We conclude with a corollary that shows how to improve the corresponding results of \cite{SSY2026_slow,SSY2026} by removing the time gap.
\begin{corollary} \label{corollary}
Let $p>1$, $q>0$  and assume $\Lambda_{p-1}>0$, where $\Lambda_{p-1}$ is defined according to \eqref{LAMBDA}. Let $u$ be a continuous, locally nonnegative weak solution
to \eqref{eq:PDE} in $\Omega_T$, with $K=K(x,y,t)$ satisfying
\eqref{k} and measurable in time. Let $(x_o,t_o)\in\Omega_T$
be such that $u_o:=u(x_o,t_o)>0$.  There exist constants $\mathfrak{s}\in(0,1)$, $\mathfrak{c}>1$ and $\mathfrak{a} \geq 1$, depending only on the data, such that if
\[
B_{16\varrho}(x_o)\times
\left[t_o-\mathfrak{a}u_o^{q+1-p}\varrho^{sp},\,
      t_o+\mathfrak{a}u_o^{q+1-p}\varrho^{sp}\right]\Subset\Omega_T,
\]
and the smallness condition on the negative tail \eqref{negativetaill} holds, then \eqref{supinf}, \eqref{integralbounds} and \eqref{integralbounds2} hold, with $u^{p-1}$ replaced by $u_+^{p-1}$ in the last two estimates.
\end{corollary}

\begin{proof}
Let $\delta$ be given by Lemma~\ref{lem3.2} with
$\texttt{c}$ satisfying \eqref{cq}. By Remark~\ref{1rem},
\eqref{negativetaill} gives \eqref{integralbounds}.
For $q \leq p-1$, \eqref{rr} gives \eqref{integralbounds2} and \eqref{supinfslow} gives     \eqref{supinf}. Note that neither the global nonnegativity of $u$ nor the independence of the kernel on $t$ is used.

For $q>p-1$, if $sp<N$, we use the local upper estimate
established in \cite[Proposition~1.4]{SSY2026}. If $sp\geq N$,
an estimate of the same form follows from
\cite[Theorem~1.3]{SZ26}, applied with $r=p-1\geq1$. Let $\tilde{\texttt{a}}\geq1$ denote the resulting constant, depending only on the data.
We apply this estimate in the cylinder
\begin{align}\label{cylinder_1}
Q^-_{\varrho,\tilde{\delta}\theta}(x_o,t_o)
=
B_\varrho(x_o)\times(t_o-\tilde{\delta}\theta,t_o),
\qquad
\theta=
\left(\frac{u_o}{4\tilde{\texttt{a}}}\right)^{q+1-p}
\varrho^{sp},
\qquad
u_o:=u(x_o,t_o),
\end{align}
with $\tilde{\delta} \in (0,1)$ to be chosen later. By an absorption argument as in \eqref{usu}, we obtain
\begin{align} \label{cpt1}
    u_o\leq
\frac{4\tilde{\texttt{a}}}
{\tilde{\delta}^{\frac{N+sp}{spq}}}
\left(
\int_{t_o-\tilde{\delta}\theta}^{t_o}
\int_{\mathbb R^N}
\frac{u_+^{p-1}(x,t)}
{(\varrho+|x-x_o|)^{N+sp}}
\, \dx\, \dt
\right)^{\frac1q},
\end{align}
that implies
\eqref{integralbounds2} for $u_+^{p-1}$.  
We introduce
\[
w(z,\tau):=
\frac{
u(x_o+\varrho z,\,
t_o-\tilde{\delta}\theta
+k_{\ast}^{q+1-p}\varrho^{sp}\tau)
}{k_{\ast}},
\]
which is a weak solution of \eqref{eq:PDE} for a kernel
satisfying \eqref{k} with the same constants. By the change of variable and \eqref{cpt1}, we get
\[
\begin{aligned}
\int_0^{\tau_o}P_w(\tau)\,\d\tau
&=\frac{1}{k_{\ast}^q}
\int_{t_o-\tilde{\delta}\theta}^{t_o}
\int_{\mathbb R^N}
\frac{u_+^{p-1}(x,t)}
{(\varrho+|x-x_o|)^{N+sp}}
\,\dx\,\dt \geq 
\frac{\tilde{\delta}^{\frac{N+sp}{sp}}}
{(4\tilde{\texttt{a}})^q}
\left(\frac{u_o}{k_{\ast}}\right)^q,
\end{aligned}
\]
where $\mathcal{P}_w, \mathcal{N}_w$ are defined according to \eqref{PN} with $u=w$. The time corresponding to $t=t_o$ is
\[
\tau_o=
\frac{\tilde{\delta}\theta}
{k_{\ast}^{q+1-p}\varrho^{sp}}.
\]
Choose $\varepsilon_{\ast}\in(0,1)$ so that
\[
2^{q+1}\varepsilon_{\ast}^{q}\leq c_1\min\{1,q\}
\]
and set
\begin{equation}\label{1}
k_{\ast}:=
\frac{\varepsilon_{\ast}}{4\tilde{\texttt{a}}}
\tilde{\delta}^{\frac{N+sp}{spq}}u_o
\leq u_o .
\end{equation}
Since
\[
\tau_o
=\tilde{\delta}\left(\frac{u_o}{4\tilde{\texttt{a}}k_{\ast}}\right)^{q+1-p}
=\frac{1}{\varepsilon_{\ast}^{q+1-p}}
\tilde{\delta}^{\frac{\Lambda_{p-1}}{spq}}
\]
and $\Lambda_{p-1}>0$, we can fix $\tilde{\delta}$ sufficiently small
that $2\tau_o\leq\delta$. Set $S:=2\tau_o$. For every
$\tau\in[\tau_o,S]$, the preceding integral estimate and
\eqref{hl} give $h(\tau)=2^q$.
On the other hand, by \eqref{negativetaill},
we obtain for every $\tau\in(0,2\tau_o)$ that
\[
\begin{aligned}
\int_0^\tau N_w(\xi)\,\d\xi
&=
\frac{1}{k_{\ast}^q}
\int_{t_o-\tilde{\delta}\theta}^{
t_o-\tilde{\delta}\theta
+k_{\ast}^{q+1-p}\varrho^{sp}\tau}
\int_{\mathbb R^N\setminus B_{4\varrho}(x_o)}
\frac{u_-^{p-1}(x,t)}
{|x-x_o|^{N+sp}}
\,\dx\,\dt\leq
\mathfrak m
\left(\frac{u_o}{k_{\ast}}\right)^q,
\end{aligned}
\]
and consequently,
\[
\ell(\tau)\leq \max \{1,q \} c_2
\mathfrak m\left(\frac{u_o}{k_{\ast}}\right)^q,
\qquad \text{for }\tau\in(0,2\tau_o),
\]
where $\ell$ is defined by \eqref{hl}. Taking $\mathfrak m$ sufficiently small, we thus
have  $\ell\leq1/\texttt{c}$ on
$(\tau_o,2\tau_o)$.
Applying Lemma~\ref{lem3.2} to $w$, with
$\texttt{c}$ as in \eqref{cq}, we get
\begin{equation*}
w(z,\tau)\geq
\left(2^q-\frac{2}{c_q}\right)^{1/q}>0
\end{equation*}
for a.e. $(z,\tau)\in
B_1\times(\tau_o,2\tau_o)$.
Scaling back, we obtain the forward lower bound
\[
u_o\leq
\gamma\inf_{B_\varrho(x_o)}u(\cdot,t)
\qquad\text{for all }t\in
\bigl(t_o,t_o+\tilde{\delta}\theta\bigr),
\]
with $\gamma$ depending only on the data. The backward upper bound follows by the same continuity and contradiction argument as in the proof of Theorem~\ref{timegapno}. Note that the required control of the negative tail at the new reference point follows from \eqref{negativetaill}, by comparing the corresponding tail integrals and choosing $\mathfrak m$ sufficiently small. Note that in the fast diffusion regime one can take
$\mathfrak{a}=1$ by choosing $\tilde{\delta}$ sufficiently small. The conclusion follows by choosing $\mathfrak{s}$ suitably small
in each diffusion regime.

\end{proof}

\section{Sharpness of the time configuration} \label{sec4}

\noindent
In this section we construct explicit positive local weak solutions of equation \eqref{eq:PDE} for a fixed kernel satisfying \eqref{k}. This example shows that a time-insensitive Harnack inequality 
cannot hold in general in the nonlocal setting for any $q>0$, 
even in the fast diffusion regime with $q$ satisfying the 
critical restriction \eqref{suq}, in contrast with the local framework. We divide the construction into a few key steps. We consider $N\geq 2$, $0<s<1$, $p>1$, $q>0$, and the domain $\Omega=B_{16}$, where we use the notation $B_r=B_r(0)$.

\vskip0.2cm
\noindent
\textit{Construction of the kernel.} For simplicity, we first consider the case $p=2$ and then show how to extend the argument to any $p>1$. We consider the following kernel
\begin{equation*}
    K(x,y)= 
    \begin{cases}
        |y|^{-N-2s}, &\quad \text{ if $x \in B_{16}$, $y \notin B_{32}$;} \\
        |x|^{-N-2s}, &\quad \text{ if $y \in B_{16}$, $x \notin B_{32}$;} \\
        |x-y|^{-N-2s}, & \quad \text{ otherwise}.
    \end{cases}
\end{equation*}
We observe that the kernel is symmetric, as $K(x,y)=K(y,x)$, and measureable. Moreover, we observe that if $x \in B_{16}$ and $y \notin B_{32}$ then
\begin{equation*}
    |x-y|\geq |y|-|x| > \frac{|y|}{2} \, , \qquad
    |x-y|\leq |x|+|y| < \frac{3}{2}|y| \, ,
\end{equation*}
this in turn, is stating that
\begin{equation}
    \label{Kprop1}
  2^{-N-2s}|x-y|^{-N-2s} \leq |y|^{-N-2s} \leq 2^{N+2s} |x-y|^{-N-2s} \, .
\end{equation}
In the same way, if $y \in B_{16}$ and $x \notin B_{32}$ we have
\begin{equation}
    \label{Kprop2}
    2^{-N-2s}|x-y|^{-N-2s} \leq |x|^{-N-2s} \leq 2^{N+2s} |x-y|^{-N-2s} \, .
\end{equation}
Therefore, putting together \eqref{Kprop1} and \eqref{Kprop2}, and remembering that in the other points of the space $K(x,y)= |x-y|^{-N-2s}$, we have for all $(x,y) \in \R^N\times\R^N$ with $x\ne y$ the bounds
\begin{equation*}
    2^{-N-2s} |x-y|^{-N-2s} \leq K(x,y) \leq 2^{N+2s} |x-y|^{-N-2s} \, .
\end{equation*}

\vskip0.2cm
\noindent
\textit{The explicit solution.} Consider $e_1=(1,0, \dots,0)$ the versor of the first axis and denote with $40e_1$ the point in space which has $40$ as first coordinate and 0 in all the others. Then, we consider the ball $B_1(40e_1)$ and the smooth cutoff function centered in this ball
\begin{equation*}
    \chi(x) = 
    \begin{cases}
        e^{-\frac{1}{1-|x-40e_1|^2}}, & \quad \text{ $x \in B_1(40e_1)$} \\
        0, & \quad \text{ $ x \notin B_1(40e_1)$}.
    \end{cases}
\end{equation*}
Afterwards, consider the quantity
\begin{equation*}
    Z = \int_{\R^N} \frac{\chi(y)}{|y|^{N+2s}} \, \dy > 0 \, ,
\end{equation*}
and the normalization of $\chi$
\begin{equation*}
    \psi(x)= \frac{\chi(x)}{Z} \, .
\end{equation*}

\noindent
We observe that $\psi \in C_0^{\infty}(\R^N)$, as it is the product of a cutoff function with a number, $\psi \geq 0$, $\text{supp} \, \psi \subset \overline{B_1(40e_1)}\subset B_{32}^c$, and 
\begin{equation*}
    \int_{\R^N} \frac{\psi(y)}{|y|^{N+2s}} \, \dy = \frac{1}{Z}\int_{\R^N} \frac{\chi(y)}{|y|^{N+2s}} \, \dy = 1 \, .
\end{equation*}
The function $\psi$ is a smooth, time-independent bump function supported far from the domain $\Omega=B_{16}$, with weighted mass equal to $1$. The role of $\psi$ is to introduce a spatially localized perturbation far from the domain $\Omega=B_{16}$, whose nonlocal contribution compensates for the time derivative of $u_k$ inside $\Omega$. This allows the solution to exhibit arbitrarily large temporal oscillations in $\Omega$ while still satisfying the equation.

Now, for every $k >0$ to be chosen, define, for all $(x,t) \in \R^N \times \R$
\begin{equation}
\label{defuk}
    \begin{aligned}
        &a_k(t)=e^{kt}, \quad b_k(t)= \frac{1}{2}\frac{\d}{\d t}(a_k^q)=\frac{qk}{2}e^{qkt} \\
        & u_k(x,t) =a_k(t) +b_k(t)\psi(x)
    \end{aligned}
\end{equation}

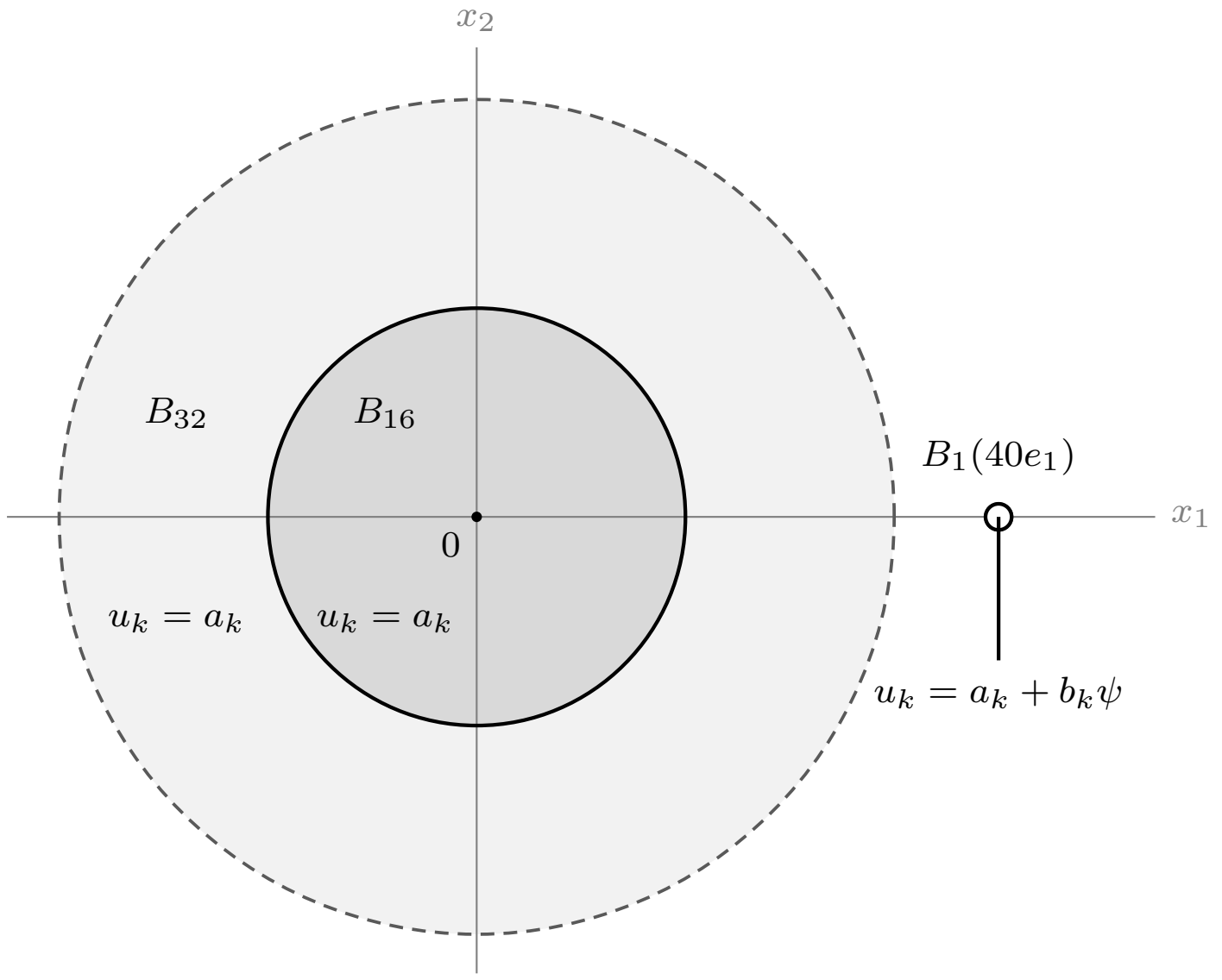
\begin{figure}[!htbp]
\centering
\begin{tikzpicture}[
    x=0.10cm, y=0.10cm,
    every node/.style={font=\footnotesize}
]
    \fill[black!5] (0,0) circle (32);
    \fill[black!15] (0,0) circle (16);

    \draw[thin, black!50] (-36,0) -- (52,0)
        node[right] {$x_1$};
    \draw[thin, black!50] (0,-35) -- (0,36)
        node[above] {$x_2$};

    \draw[black!65, line width=0.7pt,
          dash pattern=on 3pt off 2pt]
        (0,0) circle (32);

    \draw[black, line width=0.8pt]
        (0,0) circle (16);

    \fill[black!80] (40,0) circle (1);
    \draw[black, line width=0.7pt]
        (40,0) circle (1);

    \fill (0,0) circle (0.4);
    \node[below left] at (0,0) {$0$};

    \node at (-23,8) {$B_{32}$};
    \node at (-7,8) {$B_{16}$};

    \node at (-7,-8) {$u_k=a_k$};
    \node at (-23,-8) {$u_k=a_k$};

\draw[thick, fill=white] (40,0) circle (1);

\node[anchor=south] at (40,2)
    {$B_1(40e_1)$};

\node[anchor=north] (uk) at (40,-11)
    {$u_k=a_k+b_k\psi$};

\draw[thick] (40,0) -- (uk.north);

\end{tikzpicture}

\caption{Two-dimensional section of the geometry used
in the construction. The ball $B_{16}$ is the local
domain, and the bump function $\psi$ is supported in
$\overline{B_1(40e_1)}$, outside $B_{32}$.
The kernel $K(x,y)$ has a different expression
when one variable belongs to $B_{16}$ and
the other lies outside $B_{32}$.}
\label{fig:counterexample-geometry}
\end{figure}

\noindent
The globally positive solution $u_k(x,t)$ depends on space only through $\psi(x)$ and therefore coincides with $a_k$ everywhere except in $B_1(40e_1)$, while it grows exponentially in time. Therefore
\begin{equation}
\label{ulocal}
    u_k(x,t) = e^{kt} \quad \text{ in }B_{16} \quad \text{and} \quad   u_k(0,0)=1.
\end{equation}

\begin{center}

\begin{figure}[!htbp]
\includegraphics[width=0.7\linewidth]{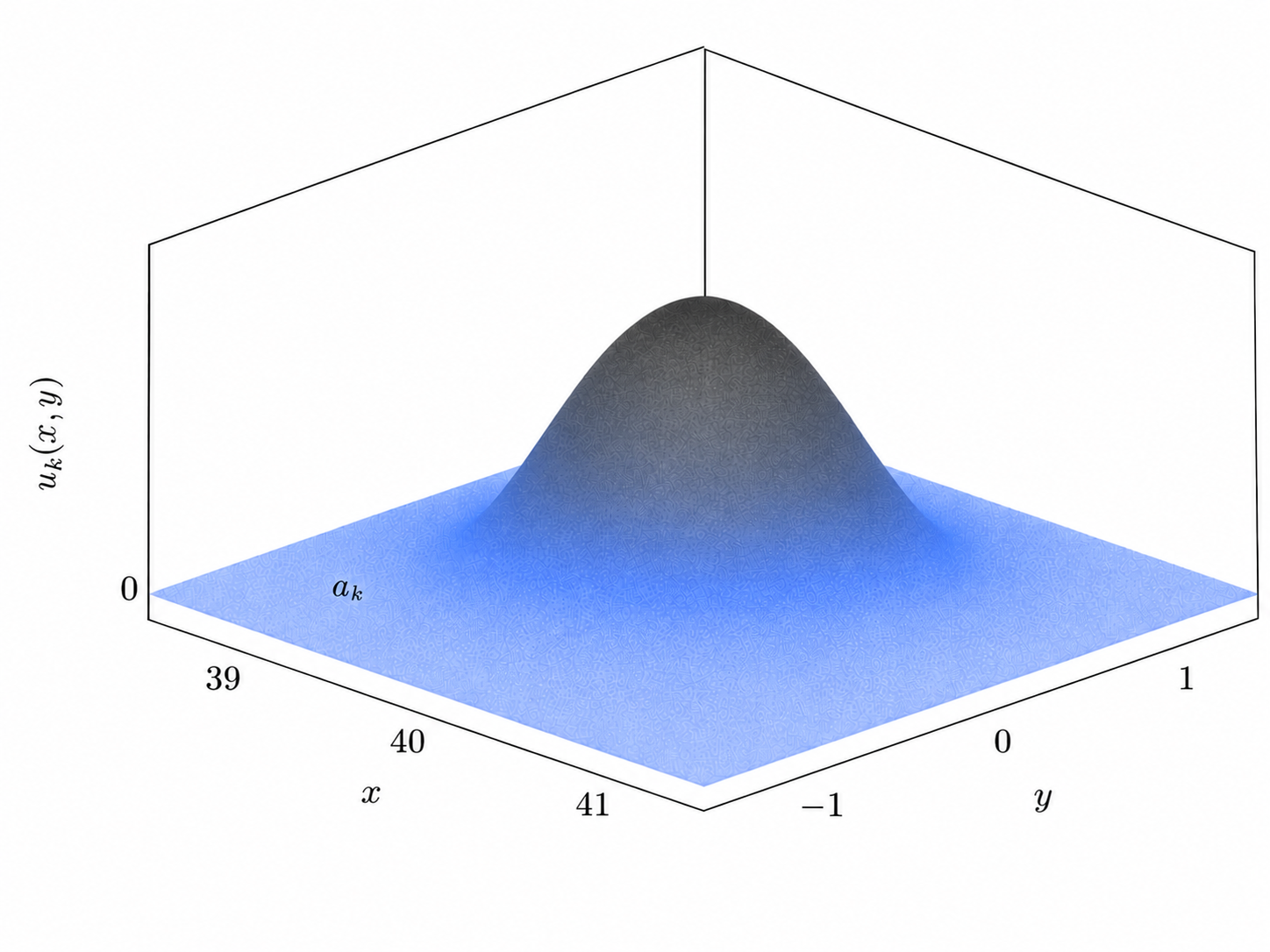}
\caption{The function $u_k$ at a fixed time $t=1$. Observe that the function attains a constant value of $a_k$ outside of $B_1(40,0)$ (the blue region), while it spikes inside of the ball (the black region). This image was made for $N=2$, $s=1/2$, $q=1$, $k=1$. The vertical scale is schematic. For these parameters the maximum is $e+1/(2Z)\approx 6.86\times 10^4$.}

\end{figure}

\end{center}

\vskip0.2cm
\noindent
\textit{Proof that $u_k$ is a local weak solution for all $k>0$.}
Now, we show that $u_k$ is indeed a local weak solution in $B_{16} \times [-T,T]$, for some $T>0$, and so it satisfies Definition \ref{def}. For the regularity properties, we observe that for every open set $D\Subset B_{16}$ and every sub-interval $[t_1,t_2]\subset [-T,T]$ we have, for every $t_o \in [t_1,t_2]$
\begin{equation*}
\begin{aligned}
    \lim_{t\to t_o} ||u_k(\cdot,t)-u_k(\cdot,t_o)||_{L^{q+1}(D)} = &\lim_{t\to t_o}  \biggl(\int_D |e^{kt}-e^{kt_o}|^{q+1} \, \dx \biggl)^{\frac{1}{q+1}} =|D|^{\frac{1}{q+1}}\lim_{t\to t_o} |e^{kt}-e^{kt_o}| = 0
    \end{aligned}
\end{equation*}
so that $u_k \in C_{\text{loc}}(-T,T;L^{q+1}_{\text{loc}}(\Omega))$. Also, it holds that, thanks to the fact that $u_k=e^{kt}$ in $D$,
\begin{equation*}
    [u_k(\cdot,t)]_{W^{s,2}(D)}=0, \qquad
    \|u_k(\cdot,t)\|_{W^{s,2}(D)}=\|u_k(\cdot,t)\|_{L^2(D)}
    =\biggl(\int_D e^{2kt} \, \dx\biggr)^{\frac{1}{2}} \, ,
\end{equation*}
and therefore
\begin{equation*}
    ||u_k||_{L^2(t_1,t_2;W^{s,2}(D))}^2 = \int_{t_1}^{t_2} \|u_k(\cdot,t)\|_{W^{s,2}(D)}^2 \, \dt = \frac{|D|}{2k}(e^{2kt_2}-e^{2kt_1})< \infty \, .
\end{equation*}
This means that $u_k \in C_{\text{loc}}(-T,T;L^{q+1}_{\text{loc}}(\Omega)) \cap L^2_{\text{loc}}(-T,T;W^{s,2}_{\text{loc}}(\Omega))$. Moreover, condition \eqref{ineq:L_1_tail_condition} holds. In fact,
\begin{equation*}
\begin{aligned}
    \int_{t_1}^{t_2} \int_{\R^N} \frac{u_k(x,t)}{1+|x|^{N+2s}}\, \dx \, \dt &=\biggl(\int_{\R^N} \frac{\dx}{1+|x|^{N+2s}}\biggl)\biggl( \int_{t_1}^{t_2}e^{kt} \, \dt \biggl) + \biggl( \int_{\R^N} \frac{\psi(x)}{1+|x|^{N+2s}} \, \dx \biggl) \biggl(\int_{t_1}^{t_2} \frac{qk}{2} e^{qkt} \, \dt\biggl) \\
  &  \leq \biggl(\int_{\R^N} \frac{\dx}{1+|x|^{N+2s}}\biggl)\biggl[\frac{e^{kt_2}-e^{kt_1}}{k} \biggl] + \biggl[\frac{e^{qkt_2}-e^{qkt_1}}{2} \biggl]<\infty \, ,
    \end{aligned}
\end{equation*}
where we have used the property of $\psi$, and the fact that
\begin{equation*}
    \int_{\R^N} \frac{\dx}{1+|x|^{N+2s}} < \infty \, .
\end{equation*}
To verify \eqref{eq}, take any $\phi \in L^2_{\text{loc}}(-T,T;W_0^{s,2}(D)) \cap W^{1,q+1}_{\text{loc}}(-T,T;L^{q+1}(D))$. It suffices to show that
\begin{equation*}
    \begin{aligned}
        -\int_{t_1}^{t_2} \int_{D} \frac{\d}{\dt}(u_k^q)(x,t)\phi(x,t) \, \dx \, \dt = \int_{t_1}^{t_2}\iint_{\R^N \times \R^N}(u_k(x,t)-u_k(y,t))(\phi(x,t)-\phi(y,t))K(x,y) \, \dx \,  \dy \, \dt.
    \end{aligned}
\end{equation*}
Taking into account the definition of $u_k$ and property \eqref{ulocal} we have that the left-hand side term becomes
\begin{equation}
\label{paraest}
\begin{aligned}
            -\int_{t_1}^{t_2} \int_{D} \frac{\d}{\dt}(u_k^q)(x,t)\phi(x,t) \, \dx \dt = -\int_{t_1}^{t_2} \int_D qke^{qkt} \, \phi(x,t) \, \dx \dt  .
\end{aligned}
\end{equation}
On the right hand-side we have divided the integral on the local and nonlocal domains
\begin{equation*}
    \begin{aligned}
        &\int_{t_1}^{t_2} \int_{B_{16}\times B_{16}} (u_k(x,t)-u_k(y,t))(\phi(x,t)-\phi(y,t))K(x,y) \, \dx\, \dy \, \dt \\
        & \qquad \qquad +2 \int_{t_1}^{t_2} \int_{B_{16}^c\times B_{16}} (u_k(x,t)-u_k(y,t))(\phi(x,t)-\phi(y,t))K(x,y) \, \dx\, \dy \, \dt \\
        & \qquad = \int_{t_1}^{t_2} \int_{B_{16}\times B_{16}} (e^{kt}-e^{kt})(\phi(x,t)-\phi(y,t))K(x,y) \, \dx\, \dy \, \dt \\
        &\qquad \qquad +2\int_{t_1}^{t_2} \int_{B_{16}^c \times B_{16}} (e^{kt} + \frac{qk}{2}e^{qkt}\psi(x)- e^{kt})(\phi(x,t)-\phi(y,t))K(x,y) \, \dx\, \dy \, \dt \\
        & \qquad = qk \int_{t_1}^{t_2} \int_{B_{16}^c \times B_{16}} e^{qkt} \psi(x) (\phi(x,t)-\phi(y,t))K(x,y) \, \dx\, \dy \, \dt \, .
    \end{aligned}
\end{equation*}
To estimate the last term, we use the fact that $\phi(\cdot,t)=0$ a.e. outside $D\subset B_{16}$ and $\operatorname{supp}\psi\subset\overline{B_1(40e_1)}\subset B_{32}^c$. Moreover, by the definition of the kernel, $K(x,y)=|x|^{-N-2s}$ whenever $y\in B_{16}$ and $x\notin B_{32}$. Hence, using the normalization of $\psi$, we obtain
\begin{equation}
\label{ellest}
    \begin{aligned}
        -qk \int_{t_1}^{t_2} \int_{B_1(40e_1)\times D}e^{qkt} \psi(x) \phi(y,t) K(x,y) \, \dx\, \dy \, \dt =&  - qk \int_{t_1}^{t_2} e^{qkt} \int_{\R^N} \frac{\psi(x)}{|x|^{N+2s}} \, \dx \int_{D} \phi(y,t) \, \dy \, \dt \\
        =& -  \int_{t_1}^{t_2} \int_{D} qke^{qkt} \phi(y,t) \, \dy \, \dt \, .
    \end{aligned}
\end{equation}
We observe that \eqref{paraest} and \eqref{ellest} coincide, therefore $u_k$ satisfies Definition \ref{def}.

Note that, for general $p>1$, replace $2s$ by $sp$ in the kernel, normalize $\int_{\R^N}\psi(y)^{p-1}|y|^{-N-sp}\, \dy=1$, and take $b_k(t)=\bigl((a_k^q)'(t)/2\bigr)^{1/(p-1)}$. Then $\mathcal L^s_{K,p}u_k=2b_k^{p-1}=(a_k^q)'$ in $B_{16}$.

\vskip0.2cm
\noindent
\textit{The Tail and the oscillation of $u_k$ in a comparison cylinder.} For simplicity, we compute the tail in the case $p=2$. The failure of time-insensitive Harnack inequalities, however, follows from the same construction for every $p>1$. We define, for some $\sigma>0 $ to be chosen, the cylinder $Q_\sigma:=B_1\times(-\sigma,\sigma]$. Let $\omega_{N-1}=|\mathbb{S}^{N-1}|$ denote the surface measure of the unit sphere. Then, due to $\text{supp} \, \psi$, we deduce
\begin{equation*}
\begin{aligned}
\text{Tail}_1(u_k;Q_\sigma) =& \int_{-\sigma}^{\sigma} \int_{B_1^c} \frac{u_k(x,t)}{|x|^{N+2s}} \, \dx \, \dt   \\
    =& \int_{-\sigma}^{\sigma} \int_{B_1^c} \frac{e^{kt}}{|x|^{N+2s}} \, \dx \dt + \frac{qk}{2} \int_{-\sigma}^{\sigma} \int_{B_1^c} \frac{e^{qkt}\psi(x)}{|x|^{N+2s}} \, \dx \dt \\
    =& \frac{\omega_{N-1}}{2s} \int_{-\sigma}^{\sigma} e^{kt} \, \dt + \frac{qk}{2} \int_{-\sigma}^{\sigma}e^{qkt} \, \dt \\
    =& \frac{\omega_{N-1}}{2sk}(e^{k\sigma}-e^{-k\sigma}) + \frac{1}{2}(e^{qk\sigma}-e^{-qk\sigma}) \,. 
    \end{aligned}
\end{equation*}
We observe that the Tail is  bounded for every $k$, but, when $k \to \infty$, it follows that
\begin{equation*}
\text{Tail}_1(u_k;Q_\sigma)^{\frac1q}
\sim
\biggl(\frac{\omega_{N-1}}{2sk}\biggr)^{\frac1q}
e^{\frac{k\sigma}{q}}
+
\biggl(\frac12\biggr)^{\frac1q}
e^{k\sigma}
\to\infty.
\end{equation*}
Moreover, we have
\begin{equation}
\label{oscbeh}
    \underset{Q_{\sigma}}{\text{ess sup}}\, u_k = e^{k\sigma}, \qquad  \underset{Q_{\sigma}}{\text{ess inf}}\, u_k = e^{-k\sigma} \, ,
\end{equation}
so that the oscillation is bounded for every $k$, but
\begin{equation*}
    \underset{Q_{\sigma}}{\text{ess osc}} \, u_k= e^{k\sigma}- e^{-k\sigma} \to \infty \, .
\end{equation*}

\noindent
\textit{Failure of time insensitive Harnack estimates.} Suppose, by contradiction, that a time-insensitive Harnack
inequality holds. More precisely, assume that there exist
constants $\mathfrak{s}\in(0,1)$, $\mathfrak{c}>1$ and
$\mathfrak{a}>1$, depending only on the data, such that every
globally positive local weak solution $u$ of \eqref{eq:PDE} satisfies
\begin{equation}\label{ellharn}
    \mathfrak{c}^{-1}u(x_o,t_o)
    \leq u(x,t)
    \leq \mathfrak{c}u(x_o,t_o)
\end{equation}
for all
\[
(x,t)\in B_\varrho(x_o)\times
\bigl(t_o-\mathfrak{s}u(x_o,t_o)^{q-1}\varrho^{2s},
      t_o+\mathfrak{s}u(x_o,t_o)^{q-1}\varrho^{2s}\bigr),
\]
whenever $\varrho>0$, $(x_o,t_o)\in\Omega'\times(0,T)$, and
\[
B_{16\varrho}(x_o)\times
\bigl[t_o-\mathfrak{a}u(x_o,t_o)^{q-1}\varrho^{2s},
      t_o+\mathfrak{a}u(x_o,t_o)^{q-1}\varrho^{2s}\bigr]
\Subset\Omega'\times(0,T).
\]

\noindent 

Now consider the local weak solution $u_k(x,t-t_o)$ defined
according to \eqref{defuk}. Choose then $\Omega'=\Omega=B_{16}$,
$\varrho=\min\{1/2,\mathfrak{a}^{-1/(2s)}\}$,
$t_o=2$, $x_o=0$ and $T=4$.
The inclusion condition is satisfied since
\begin{equation*}
    B_{16\varrho}\times
    [2-\mathfrak{a}\varrho^{2s},2+\mathfrak{a}\varrho^{2s}]
    \subset B_8\times[1,3]
    \Subset B_{16}\times(0,4),
\end{equation*}
where we used that by definition $u_k(0,0)=1$.
Now, consider the change of variables $\bar{t}=t-t_o$
and set $\sigma:=\mathfrak{s}\varrho^{2s}$.

\noindent
In this setting, if \eqref{ellharn} were to hold then,
since $u_k$ is spatially constant in $B_{16}$,
in the cylinder $Q_{\sigma}$ we would have
\begin{equation*}
    \mathfrak{c}^{-1}\leq u_k(x,t) \leq \mathfrak{c}. 
\end{equation*}
However, we observe that  
\begin{equation*}
    u_k(0,\sigma/2)= e^{k\sigma/2} \qquad
    \text{ and } \qquad
    u_k(0,-\sigma/2) = e^{-k\sigma/2}.
\end{equation*}
Therefore, choosing $k>\frac{2\ln(\mathfrak{c})}{\sigma}$,
a choice that depends only on the data, gives us
\begin{equation*}
    u_k(0,\sigma/2) > \mathfrak{c} \qquad
    \text{ and } \qquad
    u_k(0,-\sigma/2) < \mathfrak{c}^{-1}
\end{equation*}
violating both bounds. Hence, there is no possible choice of
$\mathfrak{s}$, $\mathfrak{c}$ and $\mathfrak{a}$, depending on
the data, that can make \eqref{ellharn} hold for the class
of local solutions from Definition \ref{def}.

\noindent
Moreover, we observe that, not even the classical local
time-insensitive Harnack 
\begin{equation*}
    \sup_{Q_{\sigma}} u_k
    \leq \mathfrak{c} \inf_{Q_{\sigma}} u_k
\end{equation*}
can hold, as from \eqref{oscbeh} we have
\begin{equation*}
    \frac{\sup_{Q_{\sigma}} u_k}{\inf_{Q_{\sigma}} u_k}
    = e^{2k\sigma} > \mathfrak{c}
\end{equation*}
for the choice of $k>\frac{\ln(\mathfrak{c})}{2\sigma}$
depending only on data. The calculation for general $p>1$
proceeds analogously. This completes the construction of the
counterexample, showing that time-insensitive Harnack inequalities
cannot hold in general in the nonlocal setting for any $p>1$
and $q>0$.
\vskip 0.3 cm
We observe that the failure of the time-insensitive Harnack
inequality occurs when the supremum is taken at later times
than the infimum. In fact, if we consider half time cylinders,
in the sense that the stopping point in time of the first
cylinder is precisely the starting point of the second,we have
\begin{equation*}
    Q_{\sigma}^{-}=B_1\times(-\sigma,\tau] \, ,
    \qquad Q_{\sigma}^{+}=B_1\times(\tau,\sigma]
\end{equation*}
for some $\tau\in(-\sigma,\sigma)$, as well as
\begin{equation*} 
    \sup_{Q_{\sigma}^{-}} u_k = e^{k\tau} \, ,
    \qquad \inf_{Q_{\sigma}^{+}} u_k = e^{k\tau}.
\end{equation*}
Taking $\tau=0$ in the tail calculation gives
\begin{equation*}
    \text{Tail}_1(u_k;Q_{\sigma}^{-})
    =\frac{\omega_{N-1}}{2sk}(1-e^{-k\sigma})
    +\frac12(1-e^{-qk\sigma})
    \longrightarrow \frac12 \qquad\text{as }k\to\infty,
\end{equation*}
so that, for this choice of $\tau$,
\begin{equation*}
    \frac{\sup_{Q_{\sigma}^{-}} u_k}
    {\inf_{Q_{\sigma}^{+}}u_k}
    +\frac{\text{Tail}_1(u_k;Q_{\sigma}^{-})^{1/q}}
    {\inf_{Q_{\sigma}^{+}}u_k}
    \longrightarrow 1+2^{-1/q} \, ,
\end{equation*}
which is a bounded quantity.

\vskip0.2cm
\noindent
\textit{The difference with the local case.} 
The main idea we are employing is to exploit the fractional
structure of the equations to construct solutions that could
not be considered for the local doubly nonlinear parabolic equation.
In fact, if we attempted the same construction in this case, we would have
\begin{equation*}
    \int_{t_1}^{t_2} \int_{D} \partial_t u_k^{q}
    \phi(x,t) \, \dx \dt
    = -\int_{t_1}^{t_2} \int_{D}
    |D u_k|^{p-2} D u_k \cdot D\phi(x,t) \, \dx \dt.
\end{equation*}
Since $u_k$ is spatially constant in $B_{16}$, the right-hand
side vanishes for every $D\Subset B_{16}$, whereas the left-hand
side is strictly positive for any nonnegative, nontrivial test
function. The construction of $u_k$ thus relies on the nonlocal
nature of the diffusion operator, which detects the perturbation
outside the domain where the equation is satisfied.

Observe that solutions with an exponential time dependence
are also used in the local setting to establish the optimality
of the range of exponents for which the time-insensitive Harnack
inequality holds; see \cite[Section~7.1]{B1}.

\vskip 0.3 cm
\noindent
In conclusion, the Harnack inequality without a time gap is
sharp in all the diffusion regimes with respect to the time
configuration. The solutions $u_k$ are spatially constant in
$B_{16}$, and their exponential growth rules out a uniform
comparison over the
same intrinsic cylinder.

\newpage
\noindent
\textbf{AI disclosure.} An initial draft of the counterexample was provided by ChatGPT, which was also used for editorial assistance and to assist in checking the proofs.

\end{document}